\documentclass[a4paper]{article}
\usepackage[a4paper,top=3cm,left=3cm,right=2cm,bottom=2cm]{geometry}
\usepackage[affil-it]{authblk}
\usepackage{amsmath}
\usepackage{amsthm}
\usepackage{natbib}
\usepackage{algorithm}
\usepackage{setspace}
\usepackage{fancyhdr}
\usepackage{amssymb}
\usepackage{setspace}
\usepackage[usenames,dvipsnames]{xcolor}

\usepackage{lscape}
\usepackage{subfigure}
\usepackage{graphics}
\usepackage{graphicx}
\usepackage{caption}
\usepackage{csquotes}
\usepackage[noend]{algpseudocode}
\usepackage{longtable}
\usepackage{multirow}
\usepackage{url}
\usepackage[vertfit]{breakurl} 
\usepackage{booktabs}
\usepackage{bigstrut}
\usepackage[shortlabels]{enumitem}
\usepackage{array}
\usepackage{tikz}
\usepackage{rotating}
\usepackage{morefloats}
\usepackage{colortbl}
\usepackage[normalem]{ulem}
\usepackage[titletoc,title]{appendix}

\newcommand{\DW}{Dantzig-Wolfe}
\newcommand{\PrwC}{$P|r_j|\sum w_jC_j$}
\newtheorem{property}{Property}
\usetikzlibrary{arrows}
\newcolumntype{L}[1]{>{\raggedright\let\newline\\\arraybackslash\hspace{0pt}}m{#1}}
\newcolumntype{C}[1]{>{\centering\let\newline\\\arraybackslash\hspace{0pt}}m{#1}}
\newcolumntype{R}[1]{>{\raggedleft\let\newline\\\arraybackslash\hspace{0pt}}m{#1}}



\begin{document}
\sloppy
\allowdisplaybreaks
\newcolumntype{H}{>{\setbox0=\hbox\bgroup}c<{\egroup}@{}}




\large
\title{Scheduling jobs with release dates on identical parallel machines by minimizing the total weighted completion time}

\author[1,2]{Arthur Kramer\thanks{{\footnotesize email: \texttt{arthur.kramer@ct.ufrn.br}}}}
\author[2]{Mauro Dell'Amico}
\author[3]{Dominique Feillet}
\author[2]{Manuel Iori}

\affil[1]{{\small Departamento de Engenharia de Produ\c{c}\~{a}o\\ Universidade Federal do Rio Grande do Norte, Brazil}}
\affil[2]{{\small Dipartimento di Scienze e Metodi dell'Ingegneria\\ Universit\`{a} degli Studi di Modena e Reggio Emilia, Italy}}
\affil[3]{{\small \'{E}cole des Mines de Saint-\'{E}tienne and UMR CNRS 6158 LIMOS\\ CMP Georges Charpak, France}}

\date{}

\maketitle

\vspace{-0.75cm}
\begin{center}
 Technical report -- June 2020 \\
 \footnotesize{Published version available at: {https://doi.org/10.1016/j.cor.2020.105018}}
\end{center}

\begin{abstract}
	
This paper addresses the problem of scheduling a set of jobs that are released over the time on a set of identical parallel machines, aiming at the minimization of the total weighted completion time. This problem, referred to as {\PrwC}, is of great importance in practice, because it models a variety of real-life applications. Despite its importance, the {\PrwC} has not received much attention in the recent literature. In this work, we fill this gap by proposing mixed integer linear programs and a tailored branch-and-price algorithm. Our {branch-and-price} relies on the decomposition of an arc-flow formulation and on the use of efficient exact and heuristic methods for solving the pricing subproblem. Computational experiments carried out on a set of randomly generated instances prove that the proposed methods can solve to the proven optimality instances with up to $200$ jobs and $10$ machines, and provide very low gaps for larger instances.	
\end{abstract}

%
\onehalfspace
\section{Introduction}
\label{sec:intro_PrWC}

Scheduling problems are among the most important problems in the combinatorial optimization field, being intensively studied since the early 1950s. These problems are relevant because of their practical and theoretical importance. Many practical situations are modeled by means of scheduling problems, including personnel (see, e.g., \citealt{Vandenbergh2013}), production (see, e.g., \citealt{Allahverdi2015}), health care (see, e.g., \citealt{hall2012}) and project (see, e.g., \citealt{Demeulemeester2002}) problems. In addition to their practical importance, solving scheduling problems is usually a challenging task because many of them are $\mathcal{NP}$-hard. Thus, the use of advanced optimization techniques such as column generation, branch-and-bound and/or Benders decomposition, just to cite some, is fundamental for their efficient solution.

In this work, we address a particular scheduling problem  in which we are given a set $J = \{1,\dots,n\}$ of jobs to be processed, without preemption, on a set $M = \{1,\dots,m\}$ of identical parallel machines. Each job $j \in J$ has a processing time $p_j$, a release date $r_j$ and a weight $w_j$. The objective is to find a schedule for which the sum of the weighted completion times of the jobs, $\sum_{j \in J} w_jC_j$, is a minimum, where the completion time of a job, $C_j$, is defined as the time in which its processing is finished. Using the scheduling classification of \citet{Graham1979}, this problem can be referred to as {\PrwC}. The {\PrwC} is $\mathcal{NP}$-hard because  {generalizes} the $1|r_j|\sum C_j$, which was proven {to be} $\mathcal{NP}$-hard by \cite{lenstraetal1977}.

{In the scheduling literature, only a limited number of works address the {\PrwC}. Among them, \cite{baptisteetal2008} propose  lower bounds  and \cite{nessahetal2008} a branch-and-bound algorithm. We are not aware of any recent work proposing a more efficient exact method, thus representing a lack in the literature that we try to fill with this paper.}

A variety of combinatorial optimization problems have been successfully solved by means of decomposition techniques, i.e., methods that rely on decomposing the problem in master and subproblems (see, e.g., \citealt{ValeriodeCarvalho1999, Feillet2010, Delormeetal2017} and \citealt{KowalczykLeus2018}). In this article, we focus on the application of these techniques to the {\PrwC}. In particular, we present a {\it time-indexed} (TI) formulation (see, e.g., \citealt{SousaWolsey1992}), an {\it arc-flow} (AF) formulation (see, e.g., \citealt{ValeriodeCarvalho1999}; \citealt{Delorme2016} and \citealt{Krameretal2018}), and a tailored {\it branch-and-price} (B\&P) algorithm, and use them to exactly solve {large sized} instances of the {\PrwC}. Our B\&P is based on the {\it {\DW}} (DW) decomposition of the AF formulation. On  one hand, AF formulations are characterized by having a pseudo-polynomial number of variables and constraints, and hence their application to large size instances may be prohibitive. On the other hand, the DW decomposition makes it possible to obtain a reformulation with very few constraints but with an exponential number of variables. The B\&P algorithm consists in solving at each node of a {\it branch-and-bound} (B\&B) tree a linear relaxation of the DW reformulation by means of a {\it column generation} (CG) algorithm. To this aim, we propose CG algorithms that combine the use of dominance rules and the use of exact and heuristic methods to tackle the pricing subproblem. We performed computational experiments on randomly generated instances to evaluate the performance of the proposed methods. {The results that we obtained show that our best proposed method is able to solve to the proven optimality instances with up to $200$ jobs and $10$ machines.

	The remainder of this work is organized as follows. A literature review on the {\PrwC} and closely related problems is provided in Section \ref{sec:lit_review_PrWC}. {\it Mixed integer linear programming} (MILP) TI and AF formulations are presented in Section \ref{sec:math_model_PrWC}. {The proposed B\&P algorithm is detailed in Section \ref{subsec:CG_PrwC}}. The performance evaluation of the proposed methods is presented in Section \ref{sec:comput_res_PrwC}, and then the concluding remarks and further research directions are reported in Section \ref{sec:coclusion_PrwC}.

	\section{Literature review}
	\label{sec:lit_review_PrWC}
	
	The {\PrwC} considers some important characteristics of production scheduling environments that are commonly faced in practice, such as the presence of multiple machines, release dates and priorities on the jobs. The problem has not received much attention in the recent literature, but the literature on closely related problems is rich and gives some interesting insights also on how to tackle the {\PrwC}. In the following, we briefly revise scheduling problems with release dates and involving objective functions related to the total weighted completion time minimization.
	
	The $1|r_j|\sum w_jC_j$ is a single machine $\mathcal{NP}$-hard problem for which heuristic methods, lower bounding schemes, B\&B approaches and MILP formulations have been proposed in the literature. \cite{hariripotts1983} proposed a simple heuristic method that consists in iteratively sequencing the jobs according to the non-increasing order of the $w_j/p_j$ ratio, in such a way that the release dates are satisfied and unnecessary idle times are avoided. The authors also proposed a lower bounding scheme based on the combination of {their} heuristic algorithm with a Lagrangean relaxation of the release date constraints. In addition, a B\&B integrating these methods with dominance rules was developed. With {the resulting} method, the authors solved instances with up to $50$ jobs. \cite{dyerwolsey1990} provided a theoretical analysis and developed MILP formulations to compute lower bounds for the $1|r_j|\sum w_jC_j$. Lower bounding techniques were also proposed by \cite{Belouadahetal1992}. Their main technique is based on the idea of dividing the jobs in pieces so that an efficiently solvable problem is derived.
	
	Still with regard to the $1|r_j|\sum w_jC_j$, a TI formulation was proposed by \cite{akkerEtAl2000}. The main advantage of this formulation concerns the strength of its linear relaxation, while its main drawback is due to its  pseudo-polynomial size. To work around this fact, the authors applied a DW decomposition and solved the linear relaxation of the problem by means of column and cut generation techniques, thus obtaining a valid lower bound. \cite{avellaetal2005} focused on obtaining small duality gaps. They combined a TI formulation with Lagrangean relaxation techniques to obtain near-optimal solutions for large sized instances with up to $400$ jobs. A framework featuring B\&B, {\it dynamic programming} (DP) and constraint programming techniques was proposed by \cite{panShi2008}. With this framework the authors were able to solve instances with up to $200$ jobs to the proven optimality. Later, \cite{tanakaFujikuma2012} developed a DP based algorithm that considers state elimination, dominance rules and Lagrangean relaxation concepts. To the best of our knowledge, this method represents the state-of-the-art approach for solving the $1|r_j|\sum w_jC_j$, {being able to solve to the proven optimality instances with up to $200$ jobs much faster than Pan and Shi's method}.
	
	Concerning identical parallel machine variants, many authors studied problems where preemption is allowed. Among them, \cite{Duetal1990} and \cite{Baptisteetal2007} investigated the complexity of the problem of minimizing the total flow time when preemption is allowed, where the flow time $F_j$ of a job is defined as the difference between the time in which its processing is finished and its release date, i.e., $F_j = C_j - r_j$. The former work proved that the $P|r_j, pmtn|\sum F_j$ is $\mathcal{NP}$-hard, while the latter demonstrated that if all the processing times of the jobs are identical (i.e., $p_j=p$ for all $j \in J$) the problem is solvable in polynomial time. A similar result has been achieved by \cite{BruckerKravchenko2005}, who proposed a polynomial algorithm for solving the closely related $P|r_j, p_j = p, pmtn|\sum C_j$. Regarding {again} the $P|r_j, pmtn|\sum F_j$, \cite{LeonardiRaz2007} analyzed the {\it shortest remaining processing time} (SRPT) rule, which at each time assigns to the machines the jobs already released according to their shortest remaining processing time. They showed that the SRPT rule is an $\mathcal{O}(\log(\min(n/m, p_{\max})))$-approximation algorithm for the $P|r_j, pmtn|\sum F_j$,  where {$p_{\max} = \max_{j\in J} \{p_j\}$}.
	
	Regarding parallel machine problems without preemption, \cite{LeonardiRaz2007} developed an approximation algorithm for the $P|r_j|\sum F_j$ that is based on their previous results for the problem with preemption. Their method relies on transforming a preemptive schedule into a non-preemptive one. \cite{YalaouiandChu2006}, in turn, tackled the problem of minimizing the related total completion time of the jobs. In their work, B\&B methods that consider theoretical properties of the problem as well as lower and upper bounding schemes based on job splitting and release dates relaxation are presented. The authors state that their work was the first to derive an exact approach to the $P|r_j|\sum C_j$. The results that they obtained show that instances with up to $120$ jobs and $5$ machines could be solved to the proven optimality.
	
	With regard to the weighted version, i.e., the {\PrwC} studied in this work, we highlight the contributions by \cite{Halletal1997}, \cite{baptisteetal2008} and \cite{nessahetal2008}. \cite{Halletal1997} proposed approximation algorithms with performance guarantee for a variety of scheduling problems to minimize the total weighted completion time, including the {\PrwC}. For this problem, in particular, their method is based on an approximation algorithm for the variant without release dates. The work of \cite{baptisteetal2008} focused on the study of existing and new lower bounding schemes for a variety of parallel machine scheduling problems, including the {\PrwC}. \cite{nessahetal2008}, in turn, proposed dominance rules, lower bounding schemes and {a heuristic} method. They showed that the developed lower bounds and the heuristic method provide an average gap of around $3\%$ for instances with up to $500$ jobs and $9$ machines. When these components are integrated into a B\&B framework, instances with up to $100$ jobs and $5$ machines are solved to the proven optimality. To the best of our knowledge, this method is the state-of-the-art exact algorithm for the {\PrwC}.

	\section{Mathematical formulations}
	\label{sec:math_model_PrWC}
	
	In this section, we present three different MILP formulations for {\PrwC}, namely, a time-indexed, an arc-flow and a set covering formulation.
	
	\subsection{Time-indexed formulation}
	\label{sec:TI_PrWC}
	
	{\it Time-indexed} (TI) formulations (see, e.g., \citealt{SousaWolsey1992}) seek for scheduling jobs over a sufficiently large time horizon by making use of time-indexed binary variables. Such variables, referred to as $x_{jt}$, assume value $1$ if job $j \in J$ starts its processing at time $ t \in \{0,\dots, T-p_j\}$, and $0$ otherwise, where $T$ represents the end of the time horizon. This time horizon should be sufficiently large to ensure that there exists at least an optimal solution for the problem, but needs to be as short as possible in order to reduce the formulation size. In our implementation we adopted
	\begin{align}
	T = \left\lfloor \frac{1}{m} \sum_{j \in J} p_j  + \frac{(m-1)}{m} p_{\max}\right\rfloor + r_{\max} \label{eq:time_horizon_PrWC}
	\end{align}
	{where $p_{\max}$ and $r_{\max}$  represent the maximum processing time and the maximum release date among all jobs $j \in J$, respectively. Note that Equation \eqref{eq:time_horizon_PrWC} adds $r_{\max}$ to the time horizon estimation of \cite{Krameretal2018} for the $P||\sum w_jC_j$. {Because this estimation guarantees to find an optimal solution for the $P||\sum w_jC_j$ (as also proved by \citealt{vandenAkkeretal1999}) and because all jobs are available at time $r_{max}$, we are guaranteed to find an optimal solution within $T$}.} 
	
	A TI formulation for the {\PrwC} is as follows.
	\begin{align}
	(\mbox{TI}) \quad \min \sum_{j \in J} \left( w_j p_j + \sum_{t = r_j}^{T-p_j} w_j t x_{jt} \right) \label{FO:TI_PrWC}\\
	\text{s.t.} \sum_{t = r_j}^{T - p_j} x_{jt} \geq 1 & &&  j \in J  \label{constr1:TI_PrWC}\\
	\sum_{j \in J} { \sum_{s = l_{jt}}^{u_{jt}} } x_{j s} \leq m & &&  t \in \{0,\dots,T\} \label{constr2:TI_PrWC}\\
	x_{jt} \in \{0,1\} & &&  j \in J,   t \in \{r_j, \dots, T-p_j\} \label{constr3:TI_PrWC}
	\end{align}
	where $l_{jt} = \max\{r_j, t+1-p_j\}$ and $u_{jt} = \min\{t,T+1-p_j\}$. The objective function \eqref{FO:TI_PrWC} seeks for the minimization of the total weighted completion time, constraints \eqref{constr1:TI_PrWC} ensure that all jobs are processed, constraints \eqref{constr2:TI_PrWC} guarantee that at most $m$ jobs are processed in parallel, i.e., that at most one job at a time is processed on a machine, and constraints \eqref{constr3:TI_PrWC} define the variable domain. The TI formulation \eqref{FO:TI_PrWC}-\eqref{constr3:TI_PrWC} is characterized by $\mathcal{O}(nT)$ variables and $\mathcal{O}(n+T)$ constraints.

	\subsection{Arc-flow formulation}
	\label{sec_AF_PrWC}
	
	{\it Arc-flow} (AF) formulations model problems using flows on a pseudo-polynomial size capacitated network. This technique has been successfully employed to model classical $\mathcal{NP}$-hard optimization problems, such as cutting stock problems (\citealt{ValeriodeCarvalho1999, Delorme2016}) and scheduling problems (\citealt{Krameretal2018, MradAndSouayah2018}).
	
	To formulate the {\PrwC} using an AF formulation, let us define an acyclic directed multigraph $G=(N,A)$, and a set $R$ containing all distinct job release dates plus $\{T\}$. The set of vertices $N\subseteq \{0, 1, \dots, T\}$ represents the possible starting times and completion times of a job in the planning horizon, for non dominated feasible solutions. {These times can be obtained by using Algorithm \ref{alg:N_A_PrWC}, and are usually called \textit{normal patterns} in the literature (see, e.g., \citealt{CI16}).}  The set of arcs $A$ is partitioned in such a way that $A = \cup_{j \in J} A_j \cup A_0$, where $A_j = \{(q,r,j) : q \in N, r_j \leq q \leq T-p_j \: \text{and} \: r = q + p_j\}$ is the set of job arcs associated to job $j$ and $A_0 =\{(q,r,0) : q \in N, r = \min{(t \in R: t > q)}\}$ is the set of loss arcs. Each job arc in $A_j$ represents a possible processing of a job $j$ on a machine, whereas loss arcs in $A_0$ model machine idle times. The procedure adopted to obtain {set $A$} is detailed in Algorithm \ref{alg:N_A_PrWC}.
	\begin{algorithm}[htb]
		\caption{Construction of $G=(N,A)$} \label{alg:N_A_PrWC}
		\begin{algorithmic}[1]
			\Procedure {Create\_Graph}{$T$}
			\State {\bf{initialize}} {$R\gets \text{distinct values  in } \{r_1,\dots,r_n\}\cup \{T\}$}  \text{sorted in increasing order}; 
			\State {\bf{initialize}} $i \gets 0;$
			\State {\bf{initialize}} $N \gets \emptyset; A[0 \dots n] \gets \emptyset;$ {\footnotesize {\color{gray} \Comment $N$: set of {vertices}; {$A[j]$: set $A_j$}}}
			\State {\bf{initialize}} $P[0 \dots T] \gets 0;$ {\footnotesize {\color{gray} \Comment $P$: array of size $T+1$, used to store {arc tails}}}
			\For {$q \in R$}{ $P[q] \gets 1$;} \EndFor
			\For {$q \gets 0$ to $T$}
			\If {$P[q] = 1$}
			\For {$j \gets 1$ to $n$}
			\If { {$q \geq r_j \text{ and } q + p_j \leq T$}}  { $P[q + p_j] \gets 1; A[j] \gets A[j] \cup \{(q, q + p_j,j)\}$};  \EndIf
			\EndFor
			\EndIf
			\EndFor
			\For {$q \gets 0$ to $T$}
			\If {$P[q] = 1$}
			\State {$N \gets N \cup \{q\}$};
			\While {$q \geq R[i]$} {$i \gets i+1;$}
			\EndWhile
			\State $A[0] \gets A[0] \cup \{(q,R[i],0)\};$ {\footnotesize {\color{gray} \Comment {$A[0]$}: set of dummy/loss arcs}}
			\EndIf
			\EndFor
			\State ${A \gets \cup_{j \in J} A[j]} \cup A[0]$;
			\State \bf{return} ${N,A }$
			\EndProcedure
		\end{algorithmic}
	\end{algorithm}
	
	In the AF formulation, a valid schedule of jobs on a machine corresponds to a path in {$G$}. Thus, we associate with each job arc a binary variable $x_{qrj}, (q,r,j) \in A_j$, assuming the value $1$ if job $j$ starts its processing at time $q$ and finishes at time $r$, $0$ otherwise, and with each loss arc a continuous variable $x_{qr0}, (q,r,0) \in A_0$. In this way,  the AF formulation models the {\PrwC} as the problem of finding $m$ paths in $G$ containing all jobs $j \in J$ and having a minimum total weighted completion time, namely:
	\begin{align}
	(\mbox{AF}) \quad \min \sum_{j \in J} \left( w_j p_j + \sum_{(q,r,j) \in A_j } w_j q x_{qrj} \right)\label{FO:AF_PrwC}\\
	\text{s.t.} \sum_{(q,r,j) \in A_j} x_{qrj} \geq 1 & &&  j \in J  \label{constr1:AF_PrwC}\\
	\sum_{(q,r,j)\in A} x_{qrj} - \sum_{(p,q,j)\in A} x_{pqj} = \left\{
	\begin{array}{l l}
	m, & \text{ if $q = r_{\min}$}\\
	-m, & \text{ if $q = T$}\\
	0, & \text{otherwise}
	\end{array}\right. &&& \displaystyle q \in N  \label{constr2:AF_PrwC}\\
	x_{qrj} \in \{0,1\} &&& (q,r,j) \in A \setminus A_{0} \label{constr3:AF_PrwC}\\
	0 \leq x_{qr0} \leq m &&& (q,r,0) \in A_{0} \label{constr4:AF_PrwC}
	\end{align}
	{with $r_{\min} = min_{j \in J} \{r_j$\}.}
	{The objective} function \eqref{FO:AF_PrwC} minimizes the sum of the weighted completion times; constraints \eqref{constr1:AF_PrwC} state that each job should be processed at least once; constraints \eqref{constr2:AF_PrwC} impose the flow conservation at each time instant of the time horizon; and constraints \eqref{constr3:AF_PrwC} and \eqref{constr4:AF_PrwC} define variable domains. The formulation has a pseudo-polynomial size, with $\mathcal{O}(n|N|)$ variables and $\mathcal{O}(n+|N|)$ constraints.
	
	\subsection{Set covering formulation}
	\label{subsec:SC_PrwC}
	
	{\it Set covering} (SC) formulations have been widely used to model many combinatorial optimization problems, in particular after the application of the {DW} decomposition of a given formulation. 
	Let $S$ be the set of all possible elementary schedules for a machine, i.e., those schedules where all selected jobs respect their release dates, are scheduled at most once and are not overlapping among them. Formally, $S$ is defined as the set of all solutions of polytope \eqref{constr1:AF_PrwCb}-\eqref{constr4:AF_PrwCb}:
	\begin{align}
	\sum_{(q,r,j) \in A_j} x_{qrj} \leq 1 & &&  j \in J  \label{constr1:AF_PrwCb}\\
	\sum_{(q,r,j)\in A} x_{qrj} - \sum_{(p,q,j)\in A} x_{pqj} = \left\{
	\begin{array}{l l}
	1, & \text{ if $q = r_{\min}$}\\
	-1, & \text{ if $q = T$}\\
	0, & \text{otherwise}
	\end{array}\right. &&&  q \in N \label{constr2:AF_PrwCb}\\
	x_{qrj} \in \{0,1\} &&& (q,r,j) \in A \setminus A_{0} \label{constr3:AF_PrwCb}\\
	0 \leq x_{q,r,0} \leq 1 &&& (q,r,0) \in A_{0} \label{constr4:AF_PrwCb}
	\end{align}
	
	{For a given $s \in S$, let $a_{jq}^s \in \{0,1\}$ indicate whether job $j$ starts at time $q$ in schedule $s$ or not. Then,  $a_j^s = \sum_{(q,r,j) \in A_j} a_{jq}^s$ determines if job $j$ is included in the schedule ($a_j^s =1$) or not ($a_j^s =0$), and  $c_s = \sum_{j \in J} \sum_{(q,r,j) \in A_j} w_j \: a_{jq}^s \: (q+p_j)$ is the cost of the schedule.}
	
	By defining $\pi_s, \forall s \in S$, as a binary variable assuming value $1$ if schedule $s$ is in the solution and $0$ otherwise, the {\PrwC} can be modeled as:
	\begin{align}
	(\mbox{SC}) \quad \min {\sum_{s \in S} c_s \pi_s} \label{FO:SC_PrwC}\\
	\text{s.t.} \sum_{s \in S} a_{j}^s \pi_s \geq 1 & &&  j \in J  \label{constr1:SC_PrwC}\\
	\sum_{s \in S} \pi_s \leq m & && \label{constr2:SC_PrwC}\\
	\pi_{s} \in \{0,1\} & &&  s \in S \label{constr3:SC_PrwC}
	\end{align}
	
	{The objective} function \eqref{FO:SC_PrwC} aims at minimizing the total weighted completion time; constraints \eqref{constr1:SC_PrwC} impose that each job $j \in J$ has to be processed at least once; constraint \eqref{constr2:SC_PrwC} imposes that no more than $m$ schedules are used; and constraints \eqref{constr3:SC_PrwC} define the variable domain. 
	
	As usual with SC formulations, the above model involves an exponential number of variables. In such cases, full enumeration of all variables is unpractical, and a CG approach must be used (see, e.g., \citealt{Lubbecke2005, Feillet2010}). In the following sections, we detail the B\&P method that we developed to solve this model.
	
	\section{Branch-and-price}
	\label{subsec:CG_PrwC}
	
	\subsection{Column generation principle}
	\label{subsec:CG_Principle}
	
	Column generation algorithms are often employed to solve {\it linear programs} (LPs) containing an exponential number of variables. These methods benefit from primal/dual information to avoid the full enumeration of columns.
	Basically, a CG algorithm consists of (a) solving a restricted LP containing a subset of variables; (b) retrieving dual information about the solution, identifying and adding new attractive variables (i.e., variables with a negative reduced cost, for a minimization problem) to the restricted LP. Steps (a) and (b) are repeated until no new attractive variable exists. In step (b), identifying new variables is not always a simple task. In general, this step requires solving a subproblem, known as pricing problem. This section details the components of our proposed CG method employed to solve the LP relaxation of formulation \eqref{FO:SC_PrwC}-\eqref{constr3:SC_PrwC}.
	
	First, let us define $\mbox{MP}$ as the LP relaxation of model SC (i.e., \eqref{FO:SC_PrwC}-\eqref{constr3:SC_PrwC}), referred to as master problem, and $\mbox{MP}(S')$, where $S' \subseteq S$, as the restricted master problem obtained using only the variables in $S'$. Let us also define the dual of $\mbox{MP}(S')$, named $\mbox{DMP}(S')$, as:
	\begin{align}
	(\mbox{DMP}(S')) \quad \max {\sum_{j \in J} \lambda_j + m \lambda_0}\label{FO:SubPricing_PrwC}\\
	\text{s.t.} \sum_{j \in J} a_j^s \lambda_j + \lambda_0 \leq c_s & && s \in S' \label{constr1:SubPricing_PrwC}\\
	\lambda_j \geq 0 & &&  j \in J \label{constr2:SubPricing_PrwC}\\
	\lambda_0 \leq 0 & &&  \label{constr3:SubPricing_PrwC}
	\end{align}
	where $\lambda_j, j \in J$, and $\lambda_0$ are the dual variables associated with constraints \eqref{constr1:SC_PrwC} and \eqref{constr2:SC_PrwC}, respectively.
	
	The idea of the CG algorithm is to compute the optimal solution $\pi^*$ to $\mbox{MP}(S')$ to find the corresponding dual optimal solution $\lambda^*$ and look for a column $s \in S \setminus \{S'\}$ with negative reduced cost, if any. Once there is no column with negative reduced cost, $\pi^*$ is optimal for $\mbox{MP}(S)$ and the process ends. The resulting CG procedure is shown in Algorithm \ref{alg:CG_PrwC}.
	\begin{algorithm}[htb]
		\caption{Column generation}\label{alg:CG_PrwC}
		\begin{algorithmic}[1]
			\Procedure {CG}{}
			\State {\bf{initialize}} $S' \subseteq S$ with an initial subset of schedules; \label{step:init_cols}
			\Repeat 
			\State {\bf{solve}} $\mbox{MP}(S')$ and obtain the optimal dual solution $\lambda^*$; \label{step:solve}
			\If {$\exists \: s \in S \setminus \{S'\}$ with negative reduced cost}
			$S' \gets S' \cup s$; \EndIf \label{step:add}
			\Until {$s = \emptyset$}
			\State \bf{return} $\lambda^*$
			\EndProcedure
		\end{algorithmic}
	\end{algorithm}
	
	In Algorithm \ref{alg:CG_PrwC}, the set $S'$ of initial columns is obtained by a simple heuristic procedure based on the combination of the earliest release date and the {\it weighted shortest processing time} (WSPT) rules. The method consists in iteratively sequencing the jobs on machines. {At each step, the job with lowest WSPT value among those already released is sequenced on the least loaded machine. Step \ref{step:add}, in turn, represents the process of finding and adding new columns,} which is done by solving the pricing problem that we detail in the next sections.
	
	\subsection{Pricing problem}
	\label{subsec:pricing_PrwC}
	
	In our case, the pricing problem consists in finding single machine schedules with negative reduced cost $\bar{c}_s$ satisfying \eqref{constr1:AF_PrwCb}-\eqref{constr4:AF_PrwCb}, where $\bar{c}_s$ is given by:
	\begin{align}
	\bar{c}_s = c_s - \lambda_0 -\sum_{j \in J} a_j^s \lambda_j = \left( \sum_{j \in j} \sum_{(q,r,j) \in A} a_{jq}^{s} \: \left(w_j \: (q+p_j) -\lambda_j \right) \right) - \lambda_0. \label{eq:redCost_PrwC}
	\end{align}
	
	Let us consider the acyclic directed multigraph $G=(N,A)$ introduced in Section \ref{sec_AF_PrWC} and define a cost $\bar{c}_{(q,r,j)} = w_j (q+p_j) - \lambda_j$  on every arc $(q,r,j) \in A\setminus A_0$ and a cost $\bar{c}_{(q,r,0)} = 0$ on arcs in $A_0$. A machine schedule in $S$ can then be represented by a path in $G$, starting from node $r_{min}$ and ending at node $T$. Furthermore, the cost of this path is exactly the reduced cost of the machine schedule plus $\lambda_0$. 
	
	Unfortunately, the opposite is not true: a path between  $r_{min}$ and $T$ in $G$ does not necessarily represent a schedule of $S$, because some jobs might be repeated in the path (several arcs belonging to the same set $A_j$) while set $S$ is limited to elementary schedules. For that reason, the pricing problem
	must include a further constraint, given by a non renewable resource preventing to select more than one arc from the same set $A_j$, for all $j \in J$.
	Following \cite{Feillet2004}, $n$ binary resources could be introduced to manage these constraints, at the expense of possibly prohibitive computing times. 
	
	\subsection{Pricing problem with non-elementary schedules}
	
	An alternative  emerges if we disregard constraints \eqref{constr1:AF_PrwCb} and allow jobs to be taken more than once in a schedule. Let us denote by  $S_{+}$ the set of machine schedules extended to these non-elementary schedules ($S \subseteq S_{+}$). Introducing in the master problem new variables for schedules in $ S_{+} \setminus S$  enlarges the solution space and {might consequently} weaken the lower bound provided by the model. On the other hand, it simplifies a lot the pricing problem, which becomes a standard {\it shortest path problem with resource constraints} (SPPRC) in an acyclic (topologically-ordered) graph and can be solved by dynamic programming with a complexity $\mathcal{O}(nT)$. Algorithm \ref{alg:DP_A_PrwC} reports the {\it labeling correcting} (LC) algorithm that we implemented to solve this problem. In this algorithm, $F(p)$, $p \in N$, represents the best known cost among the paths that reach node $p$. An algorithm equivalent to Algorithm \ref{alg:DP_A_PrwC} {was given} in \cite{vandenAkkeretal1999}.
	
	\begin{algorithm}[htb]
		\caption{LC$_A$} \label{alg:DP_A_PrwC}
		\begin{algorithmic}[1] 
			\Procedure {LC$_A$}{$G=(N,A)$}
			\State {\bf{initialize}} $F[r_{\min}] \gets 0$; $F[p] \gets +\infty$, $p \in N\setminus \{r_{\min}\}$;%
			\ForAll {$p \in N$ in topological order}
			\If{$F[p] < F[p+1]$} $F[p+1] \gets F[p]$;{\footnotesize {\color{gray} \Comment idle time}}
			\EndIf
			\ForAll {$(p,p+p_j,j) \in A_j, j \in J$} \label{alg:DP_A_PrwC_line4}
			\State{$q \gets p+p_j$; \quad $c \gets  w_j q - \lambda_j$};
			\If{$F[p] + c \leq F[q]$} {$F[q] \gets F[p] + c$;}
			\EndIf
			\EndFor
			\EndFor
			\State \bf{return} $F[T] - \lambda_0$
			\EndProcedure
		\end{algorithmic}
	\end{algorithm}
	
	It can be observed that the contribution of a job $j$ scheduled to start at time $q$ in the cost of a schedule is given by $w_j (q+p_j) - \lambda_j$. Note also that when {constraints \eqref{constr1:AF_PrwCb} are} relaxed there are no constraints on the number of times a job can appear in a schedule. Considering these two facts, two reduction mechanisms can be derived.
	
	First, arcs $(p,q,j) \in A_j$ can be removed from graph $G$ as soon as $\bar{c}_{(pqj)}\geq 0$. Indeed, one can always replace these arcs with zero cost paths composed of arcs in $A_0$. Equivalently, defining $y_j = \lambda_j/w_j - p_j$, arcs $(p,q,j)$ with $p \geq y_j$ can be removed; $y_j$ can be interpreted as the  starting time limit until which the job is worth being processed. In particular, if $\lambda_j \leq 0$, the job is never worth being processed: no arc from $A_j$ is kept in the graph and job $j$ is not considered in the pricing problem.
	
	The following dominance criteria can be applied to further reduce the size of the graph.
	\begin{property}
		Let us consider a pair of jobs $(i,j) \in J$, with $i \neq j$ and $\lambda_i, \lambda_j>0$. If conditions:\\*[-10mm]
		\begin{eqnarray*}
			1)&&r_i \leq r_j; \label{c:1}\\
			2)&&p_i \leq p_j; \label{c:2}\\
			3)&&w_i (r_j+p_i) - \lambda_i \leq w_j (r_j+p_j) - \lambda_j; and \label{c:3}\\
			4)&&w_i (\min(T-p_j, y_j)+p_i) - \lambda_i \leq w_j (\min(T-p_j, y_j)+p_j) - \lambda_j\label{c:4}\\*[-10mm]
		\end{eqnarray*}
		are satisfied, {then} job $i$ dominates job $j$ and set $A_j$ can be removed from the graph.
	\end{property}
	Intuitively, it means that instead of using an arc in $A_j$, it is possible and preferable to use an arc in $A_i$ (followed by one or several arcs from $A_0$ if $p_i < p_j$). This dominance relation is 
	proved as follows.
	\begin{proof}
		Let us consider a schedule including an arc $(p,p+p_j,j) \in A_j$. We show that we can schedule job $i$ instead of $j$, by improving the cost.  Because of  conditions 1) and 2),  $(p,p+p_j,j)$ can be replaced with a path starting with arc $(p,p+p_i,i)$ and possibly followed by one or several zero-cost arcs from $A_0$.  By construction and from the first reduction, $p \in \{r_j,\min(T-p_j, y_j)\}$. Then, from conditions 3)  and 4), and from the linearity of the cost function, we have $\bar{c}_{p,p+p_j,j}\geq \bar{c}_{p,p+p_i,i}$, which concludes the proof.
	\end{proof}
	
	%
	
	
	Values $y_j$ and the set of dominated jobs can easily be precomputed before starting the solution of the pricing problem. Algorithm \ref{alg:DP_A_PrwC} with the graph modified according to the two reduction mechanisms forms our first pricing method, hereafter referred to as LC$_{\text{A}}$.
	
	\subsection{Pricing with non-elementary schedules without 1-cycles}
	
	
	A relevant alternative to LC$_{\text{A}}$ is obtained by limiting the search in the pricing problem to schedules without $k$-cycles. A $k$-cycle is a sub-sequence $(j_{[i]},\dots,j_{[i+k]})$ of a schedule such that $j_{[i]} {=} j_{[i+k]}$. {We denote by  $S_{k\text{-cycle}}$ the subset of $S_{+}$ limited to non-elementary schedules without $k$-cycles. Note that $S\subseteq S_{k\text{-cycle}} \subseteq S_{+}$.}

	When $k = 1$, this problem can be solved in $\mathcal{O}(n^2T)$ time by an label correcting procedure based on the DP algorithm shown in \citet{pessoaEtAl2010}. In the next sections, we refer to this {method as} LC$_P$.
	In the following, we show how to solve this problem in $\mathcal{O}(nT)$. The proposed approach extends LC$_A$ and relies on the storage of two labels (i.e., two states) at each node $q$ of the graph, representing the two best partial schedules that finish at time $q$ with different final jobs. We introduce $F_l(q), l \in \{1,2\}$, where {$F_1(q)$} stores the best partial schedule finishing at $q$, and $F_2(q)$ represents the best schedule finishing at $q$ in which the last job is different from the last job of the first schedule. A state $L$ is represented with two values: $L.val$ is the cost of the label, $L.last$ the last job. Then, for each state $F_l(q)$ of the algorithm, {two main classes of expansions can be performed: ($i$) to $q+p_j$, $j \in J$, representing that $j$ starts in $q$ and finishes in $q+p_j$; and ($ii$) to the subsequent release date $r \in R$, for each {$j \in J: q+p_j < r$}, representing the schedule of job $j$ starting in $q$ and finishing in $q+p_j$ plus an idle time from $q+p_j$ to $r$.} Let us define $nr(p)$ as the subsequent release date of $p$, where {$nr(p) = \min \{t: t > p, t \in R \cup \{T\} \}$} (initialized as shown in Algorithm \ref{alg:init_nr}). This approach, referred to as $LC_{2l}$ from now on, is detailed in Algorithm \ref{alg:DP_2l_PrwC}.
	
	\begin{algorithm}[htb]
		\caption{Update} \label{alg:UPDATE}
		\begin{algorithmic}[1]
			\Procedure {Update($L,F,q$)}{}
			\If {$L.val < F_1[q].val$} {\footnotesize {\color{gray} \Comment $L$ will become the best label at node $q$}}
			\If {$L.last \neq F_1[q].last$} {\footnotesize {\color{gray} \Comment $F_1[q]$ will become the second label}}
			\State {$F_2[q] \gets F_1[q]$;}
			\EndIf
			\State {$F_1[q] \gets L$;}
			\ElsIf {$L.val < F_2[q].val$}
			\If {$L.last \neq F_1[q].last$} {\footnotesize {\color{gray} \Comment $L$ will become the second label}}
			\State {$F_2[q] \gets L$;}
			\EndIf
			\EndIf
			\EndProcedure
		\end{algorithmic}
	\end{algorithm}
	
	\begin{algorithm}[htb]
		\caption{Initialize vector $nr$ } \label{alg:init_nr}
		\begin{algorithmic}[1]
			\Procedure {Init\_{NR}}{$\:$}
			\State {\bf{initialize}} $R[\:];\:r \gets 0$;{\footnotesize {\color{gray} \Comment $R[\:]:$ list of distinct release dates in increasing order}}
			\State {$R \gets R \cup T$;} {\footnotesize {\color{gray} \Comment include $T$ at the end of the list}}
			\ForAll {$p \in N$ in topological order}
			\If {$p = R[r]$} {$r \gets r+1$;}\EndIf
			\State {$nr[p] \gets R[r]$;}
			\EndFor
			\EndProcedure
		\end{algorithmic}
	\end{algorithm}

	\begin{algorithm}[htb]
		\caption{LC$_{2l}$} \label{alg:DP_2l_PrwC}
		\begin{algorithmic}[1]
			\Procedure {LC$_{2l}$}{$G=(N,A)$}
			\State {\bf{initialize}} $F_l[r_{\min}].val \gets 0, l \in \{1,2\}$; $F_l[p].val \gets \infty, p \in N \setminus \{r_{\min}\}, l \in \{1,2\}$;
			\State {\bf{initialize}} $F_l[p].last \gets 0$, $p \in N$, $l \in \{1,2\};$
			\State {Init\_{NR}$(\:)$;}{\footnotesize {\color{gray} \Comment initialize list $nr[\:]$}}
			\ForAll {$p \in N$ in topological order}
			\ForAll {$(p,p+p_j,j)\in A_j, j \in J$}\label{algDP2l_line_c}
			\State {$l \gets 1;$}
			\If{$F_1[p].last = j$} {$l \gets 2$} {\footnotesize {\color{gray} \Comment expansion from best label with $last \neq j$}} \EndIf
			\State {$L.val \gets F_l[p].val + (p+p_j) \: w_j - \lambda_j; L.last \gets j$;} 	
			\State {Update($L,F,p+p_j$);} {\footnotesize {\color{gray} \Comment expansion with job $j$}}
			{\State Update($L,F,nr[p+p_j]$); {\footnotesize {\color{gray} \Comment expansion with job $j$ and idle time until next release date}} }
			\EndFor \label{algDP2l_line_d}
			\EndFor
			\State \bf{return} $F_1[T] - \lambda_0$
			\EndProcedure
		\end{algorithmic}
	\end{algorithm}
	
	By defining $Z_{MP(X)}^{*}$ as the optimal objective value of $MP(X)$, where $X$ refers to the set of columns, the following relation can be stated:
	\begin{align*}
	Z_{MP(S_+)}^{*} \leq Z_{MP(S_{1\text{-cycle}})}^{*} \leq Z_{MP(S)}^{*}
	\end{align*}
	so that the lower bound obtained by taking into consideration only elementary schedules is stronger than the lower bound obtained by the {$1$-cycle} relaxation, that in turn provides a stronger lower bound than the one obtained when all non-elementary schedules are allowed.
	Following this reasoning, and considering that for both of the aforementioned relaxations of the master problem, pricing problems are efficiently solved in pseudo-polynomial time, {the two relaxations and the three algorithms LC$_A$, LC$_P$ and LC$_{2l}$ can be all considered good options for solving the problem}.
	
	\subsection{Heuristic pricing}
	\label{subsec:heur_pricing}
	
	As the CG algorithm finishes when there is no column with negative reduced cost to be added to the master problem, one could heuristically solve the pricing and invoke the exact methods only when the heuristic methods fails in finding such columns. This combined approach is employed to speed up the CG algorithm. In this sense, we propose two greedy algorithms to solve the pricing subproblems.
	
	Our first greedy method, {denoted as $H_\text{elem}$}, aims at finding columns in $S$, i.e., elementary schedules. The method is based on LC$_A$, with the difference that, for each node $q$ of the graph, a list containing the jobs previously scheduled is stored in memory. Then, for a given node $q$, only expansions to jobs not in the list are allowed. This  heuristic, detailed in Algorithm \ref{alg:Helem_PrwC}, runs in $\mathcal{O}(n^2T)$ time and returns elementary schedules. {It} differs from Algorithm \ref{alg:DP_A_PrwC} by the use of list $L[p]$ to guarantee that only elementary schedules are obtained. This list is initialized at step \ref{alg:Helem_PrwC_l1} and updated in $\mathcal{O}(n)$ at steps \ref{alg:Helem_PrwC_l2} and \ref{alg:Helem_PrwC_l3}.
	\begin{algorithm}[htb]
		\caption{$H_{\text{elem}}$} \label{alg:Helem_PrwC}
		\begin{algorithmic}[1]
			\Procedure {H\_elem}{$G=(N,A)$}
			\State {\bf{initialize}} $L[p]  \gets \emptyset, p \in N;$  {\footnotesize {\color{gray} \Comment $L[p]$: set of all jobs in the schedule associated with $p$} } \label{alg:Helem_PrwC_l1}
			\State {\bf{initialize}} $F[r_{\min}] \gets 0$; $F[p] \gets +\infty$, $p \in N\setminus \{r_{\min}\}$;
			\ForAll {$p \in N$ in topological order}
			\If{$F[p] < F[p+1]$} $F[p+1] \gets F[p]$; $\:L[p+1] \gets L[p]$;{\footnotesize {\color{gray} \Comment idle time}} \label{alg:Helem_PrwC_l2}
			\EndIf
			\ForAll {$(p,p+p_j,j) \in A_j, j \in J \setminus L[p]$}
			\State{$q \gets p+p_j$; $\:c \gets w_j q - \lambda_j$};
			\If{$F[p] + c \leq F([q]$} {$F[q] \gets F[p] + c$; \quad $L[q] \gets L[p] \cup \{j\}$} \label{alg:Helem_PrwC_l3}
			\EndIf
			\EndFor
			\EndFor
			\State \bf{return} $F[T] - \lambda_0$
			\EndProcedure
		\end{algorithmic}
	\end{algorithm}
	
	The second proposed heuristic, shown in Algorithm \ref{alg:H1c_PrwC} and referred to as $H_\text{1-cycle}$, seeks for schedules without $1$-cycles, that is, schedules in {$S_{1\text{-cycle}}$}. This method is similar to $H_\text{elem}$, but instead of storing in memory a list of previously scheduled jobs per node, it maintains only the last scheduled job. Then, for each node $q$,  expansions are limited to jobs different from the last scheduled job.
	\begin{algorithm}[htb]
		\caption{{$H_\text{1-cycle}$}} \label{alg:H1c_PrwC}
		\begin{algorithmic}[1]
			\Procedure {{H\_1-cycle}}{$G=(N,A)$}
			\State {\bf{initialize}} $L[p]  \gets \emptyset, p \in N;$ {\footnotesize {\color{gray} \Comment $L[p]$: last job in the schedule associated with $p$}} \label{alg:H1C_PrwC_l1}
			\State {\bf{initialize}} $F[r_{\min}] \gets 0$; $F[p] \gets +\infty$, $p \in N\setminus \{r_{\min}\}$;
			\ForAll {$p \in N$ in topological order}
			\If{$F[p] < F[p+1]$} $F[p+1] \gets F[p]$; $\:L[p+1] \gets L[p]$;{\footnotesize {\color{gray} \Comment idle time}} \label{alg:H1C_PrwC_l2}
			\EndIf
			\ForAll {$(p,p+p_j,j) \in A_j, j \in J \setminus L[p]$}
			\State{$q \gets p+p_j$; $\:c \gets w_j q - \lambda_j$};
			\If{$F[p] + c \leq F([q]$} {$F[q] \gets F[p] + c$; \quad $L[q] \gets j$} \label{alg:H1C_PrwC_l3}
			\EndIf
			\EndFor
			\EndFor
			\State \bf{return} $F[T] - \lambda_0$
			\EndProcedure
		\end{algorithmic}
	\end{algorithm}

	\subsection{Column generation solution method}
	\label{subsec:CG_solMethod}
	{The methods that we detailed in the previous sections are integrated together in the CG Algorithm \ref{alg:CG_PrwC_b}.
		Firstly, we iteratively invoke a heuristic approach (step \ref{step:heur}) until it fails in finding columns with negative reduced cost, then we invoke an exact method (step \ref{step:exact}). This configuration was chosen based on preliminary experiments.}
	
	\begin{algorithm}[htb]
		\caption{Column generation - solution method}\label{alg:CG_PrwC_b}
		\begin{algorithmic}[1]
			\Procedure {CG}{}
			\State {\bf{initialize}} $S' \subseteq S$ with an initial subset of schedules; \label{step:init_cols_b}
			\State {\bf{initialize}} $heur \gets true$;
			\Repeat 
			\State {\bf{solve}} ${MP}(S')$ and obtain $\lambda^*$; \label{step:solve_b}
			\If{$heur = true$}
			\State {$s' \gets$ pricing\_heuristic();} {\footnotesize {\color{gray} \Comment returns the most negative column found by the heuristic}} \label{step:heur}
			\If {$s' = \emptyset$}{ $heur \gets false$;}
			\EndIf
			\EndIf
			\If{$heur = false$}{ $s' \gets$pricing\_exact();} {\footnotesize {\color{gray} \Comment returns the most negative column}} \label{step:exact}
			\EndIf
			\State $S' \gets S' \cup s'$;
			\Until {$s' = \emptyset$}
			\State \bf{return} $\lambda^*$
			\EndProcedure
		\end{algorithmic}
	\end{algorithm}

	\subsection{Branching strategy}
	\label{sec:BP_PrwC}
	
	Once the CG termination condition is attained, the optimal solution of $MP$ may be either integer of fractional. If it is integer, it is also optimal for SC \eqref{FO:SC_PrwC}-\eqref{constr3:SC_PrwC} and the process terminates. If instead it is fractional, branching has to be performed.
	The branching rule is an important decision to be taken in a B\&P framework. In a classical branching rule, a fractional variable $\pi_s$ is chosen and two new branches are obtained by fixing $\pi_s = 1$ and $\pi_s = 0$, {respectively. Each $\pi_s$ variable is associated with a schedule, and forbidding it to be generated during pricing} is not easy to be managed without increasing the pricing complexity. To tackle this issue, a branching rule based on the variables of AF formulation \eqref{FO:AF_PrwC}--\eqref{constr4:AF_PrwC} has been {proposed} in \cite{vandenAkkeretal1999}. We  {adopt the same} rule. 
	
	Let $\pi^*$ be an optimal $MP$ solution, with $\pi_s^*$ giving the value taken by each variable $\pi_s$. If $\pi^*$ is fractional 
	we check if it satisfies {\it Theorem 2} in \cite{vandenAkkeretal1999}, which presents a set of conditions that make it possible to transform a fractional solution into an integer feasible solution having the same cost. If the theorem is not satisfied  the process is not finished and branching is needed.
	Let $C_j(s)$ represent the completion time of job $j$ in schedule $s$. We denote by $a_j = \min\{C_j(s) | \pi^*_s > 0, a_j^s \geq 1\}$ the minimal completion time of $j$ in the schedules of $\pi^*$, and $b_j = \max\{C_j(s) | \pi^*_s > 0, a_j^s \geq 1\}$ the maximal completion time. {Let also} $\bar{J} = \{j \in J | a_j \neq b_j \}$. {Set $\bar{J}$ contains jobs that can} be selected for branching. We select the job $j \in \bar{J}$ that maximizes $(b_j - a_j)$. Then, we create two branches:
	\begin{enumerate}[label=(\roman{*})]
		\item in the first branch, we forbid arcs $(p,q,j) \in A_j$ with $q \leq \lfloor \frac{a_j+b_j}{2} \rfloor $; 
		\item in the second branch we forbid arcs $(p,q,j) \in A_j$ with $q \geq \lfloor \frac{a_j+b_j}{2} \rfloor + 1$.
	\end{enumerate}
	
	In practice, {the} procedure consists in creating new release dates and deadlines on the jobs. {Notably}, it does not change the complexity of the algorithms given in Section \ref{subsec:pricing_PrwC}.
	
	Another important decision in a B\&P algorithm concerns the choice of the branching node. In our framework, the best bound strategy is {adopted:} the node with the smallest lower bound among the open nodes is chosen for branching. The main advantage {of} this strategy is that the lower bound of the B\&B tree increases quickly. However, finding feasible solutions and closing the gap can be difficult.

	\section{Computational experiments}
	\label{sec:comput_res_PrwC}
	
	In this section, we present the computational experiments carried out to evaluate the performance of the different proposed methods for the {\PrwC}. These methods are the MILP formulations TI \eqref{FO:TI_PrWC}-\eqref{constr3:TI_PrWC} and AF \eqref{FO:AF_PrwC}-\eqref{constr4:AF_PrwC} and the B\&P methods {implemented} to solve the SC formulation \eqref{FO:SC_PrwC}-\eqref{constr3:SC_PrwC}. The B\&P methods differ {one another} by the algorithms (described in Section \ref{subsec:pricing_PrwC}) employed to solve the CG pricing subproblem.
	
	{The methods have been coded in C++ and Gurobi Optimizer $8.1$ has been adopted as MILP solver imposing it to run {on} a single thread. Our experiments were run on a computer equipped with an Intel Xeon E3-1245 processor with $3.50$GHz and $32$ GB of RAM, under Ubuntu $16.04.5$ LTS operation system.} We first detail the instance generation procedure, then present the method used to {generate} initial feasible solutions, and finally present and discuss the obtained results.

	\subsection{Benchmark instances}
	\label{subsec:bench_PrwC}
	
	To the best of our knowledge, there are {no publicly available benchmark instances from the literature regarding the {\PrwC}, so we decided to} generate a new set of instances. This set was created by generalizing the scheme {previously} proposed by \cite{nessahetal2008}. The processing times $p_j$ and weights $w_j$ of each job $j$ were generated from the uniform distributions {$U[1,p_{\max}]$ and $U[1,w_{\max}]$, respectively}. Regarding the release dates of the jobs, these were created from the uniform distribution $U[0, \frac{1}{m}(\alpha \: n \: 50.5)]$, were $n$ represents the number of jobs, $m$ the number of machines and $\alpha$ is a parameter that controls the release date tightness. Note that as the $\alpha$ parameter increases, the release dates become more dispersed, i.e., the problem becomes more constrained. 
	In our generation scheme, we adopted {the following values: $\alpha \in \{0.2, 0.6, 1.0, 1.5, 2.0, 3.0\}$, $p_{\max} \in \{100, 1000\}$, $w_{\max} \in \{10\}$, $n \in \{20, 50, 100, 125, 150, 200\}$ and $m \in \{2, 3, 5, 10\}$. For each combination of such parameters, 5 instances were created,} resulting in a total of $1440$ instances.

	\subsection{Upper bound}
	\label{subsec:UB_PrwC}
	
	In all our tests, all exact methods have been provided with an initial valid upper bound given by the {\it iterated local search} (ILS) metaheuristic proposed by \cite{KramerS2015}. With the aim on quickly obtaining a valid upper bound, in our experiments we ran the ILS algorithm with a time limit of {$n/5$} seconds or for at least one complete iteration of the inner local search phase, which is performed by a randomized variable neighborhood descent procedure.

	\subsection{Computational results}
	\label{subsec:comput_res_PrwC}
	
	In this section, we present and discuss the results obtained by the proposed solution methods. For the CG based methods, if the root node is not solved to the optimality within the time limit, then the algorithm is stopped and a valid lower bound is obtained by following the method of \cite{vanderbeckanwolsey1996}, also employed by \cite{akkerEtAl2000}. In the next subsections, we separately present the results on the instances with $p_{\max}=100$ and $p_{\max}=1000$, because this parameter has a relevant impact on the algorithmic performance.

	\subsubsection{Instances with $p_{\max}=100$}
	\label{subsec:comput_res_pmax100}
	
	First, we evaluate the exact pricing methods by running the CG algorithm with the three exact methods for solving the subproblems. In Table \ref{tab:CG_exact}, these methods are confronted by means of average execution times, {\it $t(s)$}, number of generated columns, {\it cols}, and root relaxation lower bound, {\it $lb_{lp}$}, per each group of {$30$} instances with up to {$100$} jobs and $10$ machines. {The last line provide overall average values.}
	
	\begin{table}[htb]
		\centering
		\caption{Comparison of exact methods for solving the pricing subproblem (best results in boldface) -- instances with $p_{\max}=100$}
		\scriptsize
		\setlength{\tabcolsep}{0.8mm} 
		\begin{tabular}{lrHrrrrrrrrrrrr}
			\toprule
			\multirow{2}{*}{$n$} & \multirow{2}{*}{$m$} & \multirow{2}{*}{$\alpha$} & & \multicolumn{3}{c}{LC$_A$} & & \multicolumn{3}{c}{LC$_P$} & & \multicolumn{3}{c}{LC$_{2l}$}\\
			\cmidrule{5-7}\cmidrule{9-11}\cmidrule{13-15}
			&  &  & & $lb_{lp}$ & cols & $t(s)$ & & $lb_{lp}$ & cols & $t(s)$ & & $lb_{lp}$ & cols & $t(s)$\\
			\cmidrule{1-3}\cmidrule{5-7}\cmidrule{9-11}\cmidrule{13-15}    20    & 2     &       &       & 49126.8 & \textbf{294.9} & \textbf{0.1} &       & \textbf{49137.9} & 355.3 & 0.4   &       & \textbf{49137.9} & 354.8 & 0.2 \\
			& 3     &       &       & 32017.7 & \textbf{207.0} & \textbf{0.0} &       & \textbf{32023.1} & 218.1 & 0.2   &       & \textbf{32023.1} & 218.9 & 0.1 \\
			& 5     &       &       & 22258.1 & \textbf{144.1} & \textbf{0.0} &       & \textbf{22260.9} & 144.6 & 0.1   &       & \textbf{22260.9} & 144.9 & 0.0 \\
			& 10    &       &       & 13092.9 & 78.1  & \textbf{0.0} &       & \textbf{13094.0} & 70.5  & 0.0   &       & \textbf{13094.0} & \textbf{70.2} & 0.0 \\
			\cmidrule{1-2}\cmidrule{5-7}\cmidrule{9-11}\cmidrule{13-15}    50    & 2     &       &       & 269913.4 & \textbf{2027.8} & \textbf{5.3} &       & \textbf{269987.3} & 3094.4 & 58.6  &       & \textbf{269987.3} & 3148.8 & 16.0 \\
			& 3     &       &       & 182051.8 & \textbf{1336.9} & \textbf{2.6} &       & \textbf{182085.0} & 2067.0 & 27.7  &       & \textbf{182085.0} & 2041.2 & 7.6 \\
			& 5     &       &       & 115666.4 & \textbf{1003.0} & \textbf{1.3} &       & \textbf{115675.4} & 1396.0 & 11.8  &       & \textbf{115675.4} & 1399.6 & 3.5 \\
			& 10    &       &       & 63809.8 & \textbf{618.0} & \textbf{0.5} &       & \textbf{63813.9} & 729.0 & 3.0   &       & \textbf{63813.9} & 729.5 & 0.9 \\
			\cmidrule{1-2}\cmidrule{5-7}\cmidrule{9-11}\cmidrule{13-15}    100   & 2     &       &       & 1072592.8 & \textbf{7759.6} & \textbf{135.8} &       & \textbf{1072764.3} & 14351.7 & 2385.7 &       & \textbf{1072764.3} & 14370.8 & 581.5 \\
			& 3     &       &       & 711303.4 & \textbf{5319.6} & \textbf{64.0} &       & \textbf{711370.5} & 11676.3 & 1982.4 &       & \textbf{711370.5} & 12054.5 & 893.6 \\
			& 5     &       &       & 426377.6 & \textbf{3224.3} & \textbf{23.7} &       & \textbf{426398.8} & 5562.0 & 342.9 &       & \textbf{426398.8} & 5771.0 & 81.9 \\
			& 10    &       &       & 238179.3 & \textbf{2131.6} & \textbf{10.1} &       & \textbf{238182.7} & 3179.4 & 111.3 &       & \textbf{238182.7} & 3185.9 & 26.6 \\
			\cmidrule{1-3}\cmidrule{5-7}\cmidrule{9-11}\cmidrule{13-15}   \multicolumn{3}{l}{Avg.} &       & 266365.8 & \textbf{2012.1} & \textbf{20.3} &       & \textbf{266399.5} & 3570.4 & 410.3 &       & \textbf{266399.5} & 3624.2 & 134.3 \\
			\bottomrule
		\end{tabular}%
		\label{tab:CG_exact}%
	\end{table}%
	
	From Table \ref{tab:CG_exact}, it can be noticed that the CG methods with LC$_P$ and LC$_{2l}$ always provide bounds {equal to or better than the ones obtained} with LC$_A$. This result is expected and explained by the fact that LC$_P$ and LC$_{2l}$ solve the SPPRC without one-cycles as pricing subproblem, whereas LC$_A$ solves the SPPRC itself. Thus, the solution space for the first problem is included in the solution space of the second, hence resulting in a possible stronger bound. Considering the number of columns generated, the CG with LC$_A$ {is} the method that requires less columns to prove the relaxation optimality. This results is also reflected in the execution times, that are, on average, the smallest ones among the three approaches. Regarding LC$_P$ and LC$_{2l}$, the number of columns they generate are very similar, the execution times, on the contrary, are quite different, {indicating the effectiveness of LC$_{2l}$}. The CG algorithm with LC$_{2l}$ required, on average, {$134.3$} seconds to solve the selected instances, whereas the CG with LC$_P$ needed {$410.3$} seconds. Therefore, one could select LC$_A$ as a method to solve the pricing if the objective is to quickly obtain a valid bound, and choose LC$_{2l}$ if the goal is to obtain stronger bounds.
	
	{In Table \ref{tab:CG_heur+exact}, we evaluate the impact of using heuristic approaches for solving the pricing subproblems so as to speed up the CG convergence. We report the results for the combinations of LC$_A$ and LC$_{2l}$ with $H_\text{1-cycle}$ and $H_{\text{elem}}$, resulting in four methods.} For each method and group of instances, we report the average execution times, {\it $t(s)$}, the total number of generated columns, {\it cols}, and {the average percentage number of columns generated by the exact methods, {\it $\%$cols$_e$} $= 100^*( cols_{\text{exact}}/cols)$, where $cols_{\text{exact}}$ represents the number of columns generated by the considered exact algorithm.}
	
	\begin{table}[htb]
		\centering
		\caption{Impact evaluation of combining heuristic and exact methods for solving the pricing subproblem (best results in boldface) -- instances with $p_{\max}=100$}
		\scriptsize
		\setlength{\tabcolsep}{0.8mm} 
		\begin{tabular}{lrHrrrrrrrrrrrrrrrr}
			\toprule
			\multirow{3}{*}{$n$} & \multirow{3}{*}{$m$} & \multirow{3}{*}{$\alpha$} &	& \multicolumn{7}{c}{LC$_A$} &       & \multicolumn{7}{c}{LC$_{2l}$}\\
			\cmidrule{5-11}\cmidrule{13-19}
			&	&	&	& \multicolumn{3}{c}{$H_\text{1-cycle}$} &       & \multicolumn{3}{c}{$H_{\text{elem}}$} &       & \multicolumn{3}{c}{$H_\text{1-cycle}$} &       & \multicolumn{3}{c}{$H_{\text{elem}}$} \\
			\cmidrule{5-7}\cmidrule{9-11}\cmidrule{13-15}\cmidrule{17-19}
			&  &  & & cols  & $\%$cols$_e$ & $t(s)$  &       & cols  & $\%$cols$_e$ & $t(s)$  &       & cols  & $\%$cols$_e$ & $t(s)$  &       & cols  & $\%$cols$_e$ & $t(s)$ \\
			\cmidrule{1-3}\cmidrule{5-7}\cmidrule{9-11}\cmidrule{13-15}\cmidrule{17-19}    20    & 2     & 0.0   &       & 323.2 & 11.0  & 0.1   &       & 161.7 & 20.1  & 0.1   &       & 292.9 & \textbf{1.6} & 0.1   &       & \textbf{138.2} & 8.2   & \textbf{0.1} \\
			& 3     & 0.0   &       & 203.7 & 6.2   & 0.0   &       & 102.9 & 10.2  & 0.0   &       & 193.5 & \textbf{1.4} & 0.0   &       & \textbf{94.2} & 3.4   & \textbf{0.0} \\
			& 5     & 0.0   &       & 132.3 & 2.9   & 0.0   &       & 73.4  & 3.8   & 0.0   &       & 129.5 & \textbf{0.8} & 0.0   &       & \textbf{71.7} & 1.9   & \textbf{0.0} \\
			& 10    & 0.0   &       & 67.7  & 1.8   & 0.0   &       & 41.1  & 3.4   & \textbf{0.0} &       & 67.4  & \textbf{1.5} & 0.0   &       & \textbf{40.6} & 2.6   & 0.0 \\
			\cmidrule{1-3}\cmidrule{5-7}\cmidrule{9-11}\cmidrule{13-15}\cmidrule{17-19}    50    & 2     & 0.0   &       & 3163.0 & 13.7  & 11.0  &       & 1318.7 & 31.6  & 3.1   &       & 2885.0 & \textbf{4.3} & 10.2  &       & \textbf{1061.8} & 15.9  & \textbf{2.7} \\
			& 3     & 0.0   &       & 1974.5 & 6.9   & 5.1   &       & 651.0 & 14.6  & 1.2   &       & 1892.3 & \textbf{2.1} & 5.0   &       & \textbf{583.4} & 5.6   & \textbf{1.1} \\
			& 5     & 0.0   &       & 1230.8 & 3.5   & 2.1   &       & 445.3 & 6.2   & 0.6   &       & 1202.8 & \textbf{0.4} & 2.0   &       & \textbf{415.3} & 0.8   & \textbf{0.6} \\
			& 10    & 0.0   &       & 641.3 & 1.7   & 0.5   &       & 338.7 & 1.8   & 0.3   &       & 633.8 & \textbf{0.2} & 0.5   &       & \textbf{333.0} & 0.3   & \textbf{0.3} \\
			\cmidrule{1-3}\cmidrule{5-7}\cmidrule{9-11}\cmidrule{13-15}\cmidrule{17-19}    100   & 2     & 0.0   &       & 13947.7 & 17.5  & 354.1 &       & 6244.0 & 36.8  & 92.1  &       & 12787.2 & \textbf{10.2} & 321.1 &       & \textbf{5334.9} & 27.0  & \textbf{78.9} \\
			& 3     & 0.0   &       & 9510.9 & 9.5   & 167.3 &       & 3064.0 & 22.6  & 24.3  &       & 8975.6 & \textbf{3.5} & 160.0 &       & \textbf{2624.4} & 10.6  & \textbf{21.6} \\
			& 5     & 0.0   &       & 5388.9 & 3.9   & 59.4  &       & 1564.7 & 8.2   & 7.9   &       & 5296.1 & \textbf{1.8} & 58.8  &       & \textbf{1476.9} & 3.7   & \textbf{7.6} \\
			& 10    & 0.0   &       & 2884.7 & 1.0   & 18.3  &       & 1201.6 & 1.5   & 4.3   &       & 2865.9 & \textbf{0.2} & 18.2  &       & \textbf{1185.8} & 0.4   & \textbf{4.2} \\
			\cmidrule{1-3}\cmidrule{5-7}\cmidrule{9-11}\cmidrule{13-15}\cmidrule{17-19}\multicolumn{3}{l}{Avg.} &       & 3289.1 & 6.6   & 51.5  &       & 1267.3 & 13.4  & 11.1  &       & 3101.8 & \textbf{2.3} & 48.0  &       & \textbf{1113.4} & 6.7   & \textbf{9.8} \\
			\bottomrule
		\end{tabular}%
		\label{tab:CG_heur+exact}%
	\end{table}%
	
	The results shown in Table \ref{tab:CG_heur+exact} indicate that the performance of $H_{\text{elem}}$ is, on average, better than the one by $H_\text{1-cycle}$, {because of the reduced number of columns generated as well as the smaller average execution times.} For example, the CG algorithm required, on average, {$3101.8$ columns and $48$ seconds} to solve the problem when using LC$_{2l}$ and $H_\text{1-cycle}$. By keeping LC$_{2l}$ but replacing $H_\text{1-cycle}$ by $H_{\text{elem}}$, the average number of generated columns is reduced by {more than $60\%$ and the average execution time drops to only $9.8$ seconds}. In other words, $H_{\text{elem}}$ seems to be able to found ``better" columns than $H_\text{1-cycle}$, thus speeding up the convergence at the beginning of the CG algorithm. {From Table \ref{tab:CG_heur+exact}, it can also be observed that most of the columns are generated by the heuristic methods.}
	
	By {comparing} the results in Table \ref{tab:CG_heur+exact} with the ones in Table \ref{tab:CG_exact}, it can be noticed the effectiveness of the combination of exact and heuristic methods in the CG pricing problem. Taking as example the results obtained by the CG with LC$_{2l}$ in Table \ref{tab:CG_exact}, and the results by the CG with LC$_{2l}$ and $H_{\text{elem}}$ in Table \ref{tab:CG_heur+exact}, it can be seen that the average number of generated columns drops from {$3624.2$ to only $1113.4$. The same occurs with the average execution times that passes from $134.3$ to $9.8$ seconds}.
	
	{Considering the above results, we thus adopted the combined use of LC$_{2l}$ and $H_{\text{elem}}$ as the method for solving the pricing subproblem in our CG algorithm, and embedded them into our B\&P framework, as described in Section \ref{sec:BP_PrwC}. In the following, the resulting method is named B\&P$_h$.}
	In Tables \ref{tab:BP_MILP_20-100} and \ref{tab:BP_MILP_125-200}, B\&P$_h$ is compared with the TI \eqref{FO:TI_PrWC}-\eqref{constr3:TI_PrWC} and AF \eqref{FO:AF_PrwC}-\eqref{constr4:AF_PrwC} formulations in terms of number of instances optimally solved, {\it $opt$}, percentage gap for the continuous and for the best lower bound, 
	$gap_{lp}(\%)$ and $gap(\%)$, respectively, number of open nodes, $nd$, and execution times, $t(s)$. {The} results are aggregated by $n$ and $m$, where for each combination of $(n,m)$, six instances (one for each value of $\alpha \in \{0.2,0.6,1.0,1.5,2.0,3.0\}$) have been considered. For each run of TI, AF and B\&P$_h$, a time limit of $1800$ seconds has been imposed. {The last line gives the sum of the opt values, as well as average values for the other columns.}
	
	\begin{table}[htb]
		\centering
		\caption{Comparison of B\&P$_h$ performance with MILP models TI and AF - instances with $p_{\max}=100$ and $n \in \{20,50,100\}$ (best results in boldface)}
		\scriptsize
		\setlength{\tabcolsep}{0.8mm} 
		\begin{tabular}{lrrrrrrrrrrrrrrrrrrr}
			\toprule
			\multirow{3}{*}{$n$} & \multirow{3}{*}{$m$} & &  \multicolumn{5}{c}{TI} & & \multicolumn{5}{c}{AF} & & \multicolumn{5}{c}{B\&P$_h$} \\
			\cmidrule{4-8}\cmidrule{10-14}\cmidrule{16-20}&  &  & \multirow{2}{*}{$opt$} & {$gap_{lp}$} & {$gap$} & \multirow{2}{*}{$nd$}  & \multirow{2}{*}{$t(s)$} &       & \multirow{2}{*}{$opt$} & {$gap_{lp}$} & {$gap$} & \multirow{2}{*}{$nd$}  & \multirow{2}{*}{$t(s)$} &       & \multirow{2}{*}{$opt$} & {$gap_{lp}$} & {$gap$} & \multirow{2}{*}{$nd$}  & \multirow{2}{*}{$t(s)$} \\
			&  &  &  & {$(\%)$} & {$(\%)$} &  &  &       &   & {$(\%)$} & {$(\%)$} &  &  &       &   & {$(\%)$} & {$(\%)$} &  & \\
			\cmidrule{1-2}\cmidrule{4-8}\cmidrule{10-14}\cmidrule{16-20}    {20} & 2     &       & \textbf{6} & 0.083 & \textbf{0.000} & \textbf{0.2} & 1.7   &       & \textbf{6} & 0.083 & \textbf{0.000} & \textbf{0.2} & 0.3   &       & \textbf{6} & \textbf{0.055} & \textbf{0.000} & 0.7   & \textbf{$<$0.1} \\
			& 3     &       & \textbf{6} & 0.033 & \textbf{0.000} & 0.3   & 1.2   &       & \textbf{6} & 0.033 & \textbf{0.000} & \textbf{0.2} & 0.2   &       & \textbf{6} & \textbf{0.023} & \textbf{0.000} & 1.0   & \textbf{$<$0.1} \\
			& 5     &       & \textbf{6} & 0.068 & \textbf{0.000} & \textbf{0.3} & 0.6   &       & \textbf{6} & 0.068 & \textbf{0.000} & \textbf{0.3} & 0.1   &       & \textbf{6} & \textbf{0.054} & \textbf{0.000} & 1.3   & \textbf{$<$0.1} \\
			& 10    &       & \textbf{6} & \textbf{0.000} & \textbf{0.000} & \textbf{0.0} & 0.2   &       & \textbf{6} & \textbf{0.000} & \textbf{0.000} & \textbf{0.0} & $<$0.1  &       & \textbf{6} & \textbf{0.000} & \textbf{0.000} & \textbf{0.0} & \textbf{$<$0.1} \\
			\cmidrule{1-2}\cmidrule{4-8}\cmidrule{10-14}\cmidrule{16-20}    {50} & 2     &       & \textbf{6} & 0.092 & \textbf{0.000} & 938.0 & 48.3  &       & \textbf{6} & 0.092 & \textbf{0.000} & 573.2 & \textbf{19.6} &       & 4     & \textbf{0.059} & 0.007 & \textbf{80.8} & 610.3 \\
			& 3     &       & \textbf{6} & 0.083 & \textbf{0.000} & 662.2 & 24.8  &       & \textbf{6} & 0.083 & \textbf{0.000} & 859.2 & \textbf{19.4} &       & \textbf{6} & \textbf{0.054} & \textbf{0.000} & \textbf{46.5} & 44.4 \\
			& 5     &       & \textbf{6} & 0.037 & \textbf{0.000} & 221.5 & 9.3   &       & \textbf{6} & 0.037 & \textbf{0.000} & 267.0 & \textbf{3.6} &       & \textbf{6} & \textbf{0.026} & \textbf{0.000} & \textbf{58.3} & 11.3 \\
			& 10    &       & \textbf{6} & 0.007 & \textbf{0.000} & 2.2   & 2.6   &       & \textbf{6} & 0.007 & \textbf{0.000} & \textbf{1.7} & \textbf{0.4} &       & \textbf{6} & \textbf{0.003} & \textbf{0.000} & 9.5   & \textbf{0.4} \\
			\cmidrule{1-2}\cmidrule{4-8}\cmidrule{10-14}\cmidrule{16-20}    {100} & 2     &       & \textbf{4} & 0.314 & 0.256 & 3396.3 & \textbf{863.6} &       & 3     & 0.314 & 0.256 & 7385.8 & 910.8 &       & 2     & \textbf{0.264} & \textbf{0.215} & \textbf{22.3} & 1226.8 \\
			& 3     &       & \textbf{4} & 0.081 & 0.055 & 1167.0 & 734.6 &       & \textbf{4} & 0.081 & 0.055 & 2675.5 & \textbf{710.5} &       & 3     & \textbf{0.065} & \textbf{0.035} & \textbf{50.7} & 931.3 \\
			& 5     &       & \textbf{5} & 0.042 & 0.005 & 7010.0 & \textbf{383.9} &       & \textbf{5} & 0.042 & \textbf{0.002} & 21661.3 & 450.1 &       & 3     & \textbf{0.037} & 0.005 & \textbf{199.3} & 907.4 \\
			& 10    &       & \textbf{6} & 0.013 & \textbf{0.000} & \textbf{239.8} & 22.7  &       & \textbf{6} & 0.013 & \textbf{0.000} & 252.2 & \textbf{9.0} &       & 4     & \textbf{0.009} & 0.002 & 611.0 & 607.1 \\
			\cmidrule{1-2}\cmidrule{4-8}\cmidrule{10-14}\cmidrule{16-20}\multicolumn{2}{l}{Sum/Avg.} &       & \textbf{67} & 0.071 & 0.026 & 1136.5 & 174.5 &       & 66    & 0.071 & 0.026 & 2806.4 & \textbf{177.0} &       & 58    & \textbf{0.054} & \textbf{0.022} & \textbf{90.1} & 361.6 \\
			\bottomrule		
		\end{tabular}%
		\label{tab:BP_MILP_20-100}%
	\end{table}
	
	{With regard to the results obtained for the instances with up to $100$ jobs, in Table \ref{tab:BP_MILP_20-100}, it can be noticed that TI and AF perform better than B\&P$_h$ in terms of number of instances solved to the proven optimality and execution times. The B\&P$_h$ algorithm, however, provides stronger continuous relaxation bounds and requires the exploration of a smaller amount of nodes. The results for the large-sized instances with $n>100$ are given in Table \ref{tab:BP_MILP_125-200}. For the instances with $200$ jobs and $2$ machines, TI and AF suffer from memory limit, indicated by the entry ``{\it m.lim}'', due to their large (pseudo-polynomial) number of variables and constraints. We thus provide two overall Sum/Avg. lines, one referring to all instances and the other only to those for which no memory limit was encountered by a method. Summarizing, on one side B\&P$_h$ is an interesting solution method because it does not suffer from memory limits and provides very low gaps. On the other side, TI and AF are very efficient and find a higher number of proven optimal solutions overall.}

	\begin{table}[htb]
		\centering
		\caption{Comparison of B\&P$_h$ performance with MILP models TI and AF- instances with $p_{\max}=100$ and $n \in \{125,150,200\}$ (best results in boldface)}
		\scriptsize
		\setlength{\tabcolsep}{0.8mm} 
		\begin{tabular}{lrrrrrrrrrrrrrrrrrrr}
			\toprule
			\multirow{3}{*}{$n$} & \multirow{3}{*}{$m$} & &  \multicolumn{5}{c}{TI} & & \multicolumn{5}{c}{AF} & & \multicolumn{5}{c}{B\&P$_h$} \\
			\cmidrule{4-8}\cmidrule{10-14}\cmidrule{16-20}&  &  & \multirow{2}{*}{$opt$} & {$gap_{lp}$} & {$gap$} & \multirow{2}{*}{$nd$}  & \multirow{2}{*}{$t(s)$} &       & \multirow{2}{*}{$opt$} & {$gap_{lp}$} & {$gap$} & \multirow{2}{*}{$nd$}  & \multirow{2}{*}{$t(s)$} &       & \multirow{2}{*}{$opt$} & {$gap_{lp}$} & {$gap$} & \multirow{2}{*}{$nd$}  & \multirow{2}{*}{$t(s)$} \\
			&  &  &  & {$(\%)$} & {$(\%)$} &  &  &       &   & {$(\%)$} & {$(\%)$} &  &  &       &   & {$(\%)$} & {$(\%)$} &  & \\
			\cmidrule{1-2}\cmidrule{4-8}\cmidrule{10-14}\cmidrule{16-20}    {125} & 2     &       & \textbf{4} & 0.361 & \textbf{0.338} & 2852.0 & \textbf{774.0} &       & \textbf{4} & 0.361 & 0.341 & 2831.0 & 787.7 &       & 2     & \textbf{0.343} & 0.339 & \textbf{10.3} & 1303.1 \\
			& 3     &       & \textbf{4} & 0.192 & 0.175 & 83.8  & 642.2 &       & \textbf{4} & 0.192 & 0.178 & 597.5 & \textbf{612.8} &       & 3     & \textbf{0.176} & \textbf{0.169} & \textbf{35.3} & 943.8 \\
			& 5     &       & \textbf{4} & 0.084 & 0.054 & 1846.2 & \textbf{878.1} &       & 3     & 0.084 & 0.056 & 2397.0 & 902.7 &       & 3     & \textbf{0.070} & \textbf{0.043} & \textbf{141.2} & 913.2 \\
			& 10    &       & \textbf{6} & 0.015 & \textbf{0.000} & 2638.8 & 138.4 &       & \textbf{6} & 0.015 & \textbf{0.000} & 3265.2 & \textbf{91.9} &       & 3     & \textbf{0.012} & 0.002 & \textbf{700.0} & 905.4 \\
			\cmidrule{1-2}\cmidrule{4-8}\cmidrule{10-14}\cmidrule{16-20}    {150} & 2     &       & \textbf{3} & 0.287 & 0.278 & \textbf{2.0} & 965.0 &       & \textbf{3} & 0.287 & 0.281 & 277.5 & \textbf{910.0} &       & 1     & \textbf{0.278} & \textbf{0.277} & 3.7   & 1568.9 \\
			& 3     &       & \textbf{3} & 0.436 & 0.428 & 21.3  & 942.1 &       & \textbf{3} & 0.436 & 0.426 & 577.5 & \textbf{906.0} &       & \textbf{3} & \textbf{0.425} & \textbf{0.385} & \textbf{20.0} & 1329.1 \\
			& 5     &       & \textbf{3} & 0.280 & 0.274 & 854.5 & 919.9 &       & \textbf{3} & 0.280 & 0.273 & 2561.3 & \textbf{902.6} &       & \textbf{3} & \textbf{0.276} & \textbf{0.268} & \textbf{73.8} & 924.4 \\
			& 10    &       & \textbf{6} & 0.011 & \textbf{0.000} & 586.5 & \textbf{260.5} &       & \textbf{6} & 0.011 & \textbf{0.000} & 1244.2 & 285.8 &       & 4     & \textbf{0.009} & 0.003 & \textbf{307.8} & 807.1 \\
			\cmidrule{1-2}\cmidrule{4-8}\cmidrule{10-14}\cmidrule{16-20}    {200} & 2     &       & 0     &       &       &       & m.lim &       & 0     &       &       &       & m.lim &       & 0     & \textbf{0.290} & \textbf{0.290} & \textbf{0.5} & t.lim \\
			& 3     &       & \textbf{3} & 0.242 & 0.236 & 229.7 & 988.4 &       & \textbf{3} & 0.242 & 0.238 & 1505.8 & \textbf{932.4} &       & 2     & \textbf{0.234} & \textbf{0.230} & \textbf{3.0} & 1377.5 \\
			& 5     &       & \textbf{3} & 0.243 & 0.238 & \textbf{0.7} & 941.1 &       & \textbf{3} & 0.243 & 0.240 & 752.0 & \textbf{905.2} &       & \textbf{3} & \textbf{0.238} & \textbf{0.234} & 25.5  & 1073.2 \\
			& 10    &       & \textbf{4} & 0.047 & 0.043 & 470.5 & 639.0 &       & \textbf{4} & 0.047 & 0.042 & 1424.8 & \textbf{608.0} &       & 3     & \textbf{0.045} & \textbf{0.040} & \textbf{116.3} & 947.6 \\
			\cmidrule{1-2}\cmidrule{4-8}\cmidrule{10-14}\cmidrule{16-20}
			\multicolumn{2}{l}{Sum/Avg.$^*$} &       & \textbf{43} & 0.200 & 0.188 & 871.5 & 735.3 &       & 42    & 0.200 & 0.189 & 1584.9 & \textbf{713.2} &       & 30    & \textbf{0.191} & \textbf{0.181} & \textbf{130.6} & 1099.4 \\
			\multicolumn{2}{l}{Sum/Avg.} &       &  &  &  &  &  &       &     &  &  & &  &       & \textbf{30}    & \textbf{0.200} & \textbf{0.190} & \textbf{119.8} & \textbf{1157.8} \\
			\bottomrule
			\multicolumn{20}{l}{$^*$Sum and average values for all instances, but the ones with $(n,m) \in \{(200,2)\}$}
		\end{tabular}%
		\label{tab:BP_MILP_125-200}%
	\end{table}
	
	\subsubsection{{Instances with $p_{\max}=1000$}}
	\label{subsec:comput_res_pmax1000}
	
	In this section, we replicate the experiments described in Section \ref{subsec:comput_res_pmax100}, but now considering the instances with $p_{\max}= 1000$. In Table \ref{tab:CG_exact_pmax1000}, we evaluate the performance of the three exact algorithms employed to solve the pricing subproblem. The results we obtained are aligned with the ones reported in Table \ref{tab:CG_exact} for $p_{\max}=100$, thus confirming the performance gain achieved when replacing LC$_P$ by LC$_{2l}$.
	
	\begin{table}[htbp]
		\centering
		\caption{Comparison of exact methods for solving the pricing subproblem (best results in boldface) -- instances with $p_{\max}=1000$}
		\scriptsize
		\setlength{\tabcolsep}{0.8mm} 
		\begin{tabular}{lrHrrrrrrrrrrrr}
			\toprule
			\multirow{2}{*}{$n$} & \multirow{2}{*}{$m$} & \multirow{2}{*}{$\alpha$} & & \multicolumn{3}{c}{LC$_A$} & & \multicolumn{3}{c}{LC$_P$} & & \multicolumn{3}{c}{LC$_{2l}$}\\
			\cmidrule{5-7}\cmidrule{9-11}\cmidrule{13-15}
			&  &  & & $lb_{lp}$ & cols & $t(s)$ & & $lb_{lp}$ & cols & $t(s)$ & & $lb_{lp}$ & cols & $t(s)$\\	
			\cmidrule{1-3}\cmidrule{5-7}\cmidrule{9-11}\cmidrule{13-15}
			{20} & 2     &       &       & 195006.7 & \textbf{229.1} & \textbf{0.1} &       & \textbf{195141.0} & 232.0 & 1.5   &       & \textbf{195141.0} & 232.5 & 0.5 \\
			& 3     &       &       & 148947.7 & \textbf{171.1} & \textbf{0.1} &       & \textbf{149007.6} & 177.1 & 0.8   &       & \textbf{149007.6} & 176.4 & 0.2 \\
			& 5     &       &       & 102597.2 & 123.9 & \textbf{0.0} &       & \textbf{102599.0} & 121.4 & 0.3   &       & \textbf{102599.0} & \textbf{120.9} & 0.1 \\
			& 10    &       &       & 71462.7 & 76.0  & \textbf{0.0} &       & \textbf{71463.8} & 74.3  & 0.1   &       & \textbf{71463.8} & \textbf{73.4} & 0.0 \\
			\cmidrule{1-3}\cmidrule{5-7}\cmidrule{9-11}\cmidrule{13-15}
			{50} & 2     &       &       & 1049597.4 & \textbf{1425.1} & \textbf{4.8} &       & \textbf{1050223.8} & 1648.2 & 207.8 &       & \textbf{1050223.8} & 1637.4 & 18.3 \\
			& 3     &       &       & 799434.9 & \textbf{952.1} & \textbf{2.3} &       & \textbf{799746.8} & 1075.5 & 91.2  &       & \textbf{799746.8} & 1084.7 & 8.5 \\
			& 5     &       &       & 523452.1 & \textbf{647.4} & \textbf{1.0} &       & \textbf{523544.7} & 693.8 & 34.7  &       & \textbf{523544.7} & 700.9 & 3.4 \\
			& 10    &       &       & 285946.9 & \textbf{399.6} & \textbf{0.4} &       & \textbf{285969.7} & 422.8 & 11.0  &       & \textbf{285969.7} & 420.5 & 1.2 \\
			\cmidrule{1-3}\cmidrule{5-7}\cmidrule{9-11}\cmidrule{13-15}
			{100} & 2     &       &       & 4173231.9 & \textbf{6914.0} & \textbf{136.0} &       & \textbf{4174522.1} & 8231.0 & 8930.3 &       & \textbf{4174522.1} & 8401.6 & 432.7 \\
			& 3     &       &       & 2842626.6 & \textbf{3797.9} & \textbf{41.0} &       & \textbf{2843104.3} & 4421.7 & 3127.0 &       &\textbf{ 2843104.3} & 4369.0 & 135.8 \\
			& 5     &       &       & 1765972.9 & \textbf{2170.3} & \textbf{15.1} &       & \textbf{1766155.7} & 2552.4 & 1166.9 &       & \textbf{1766155.7} & 2531.9 & 52.8 \\
			& 10    &       &       & 994029.4 & \textbf{1297.6} & \textbf{4.9} &       & \textbf{994066.5} & 1435.1 & 353.8 &       & \textbf{994066.5} & 1430.9 & 16.3 \\
			\cmidrule{1-3}\cmidrule{5-7}\cmidrule{9-11}\cmidrule{13-15}
			\multicolumn{3}{l}{Avg.} &       & 1079358.9 & \textbf{1517.0} & \textbf{17.1} &       & \textbf{1079628.7} & 1757.1 & 1160.4 &       & \textbf{1079628.7} & 1765.0 & 55.8 \\	
			\bottomrule
		\end{tabular}%
		\label{tab:CG_exact_pmax1000}%
	\end{table}%
	
	In Table \ref{tab:CG_heur+exact_pmax1000}, we {evaluate} the combined use of heuristic and exact methods in the CG. 
	\begin{table}[htbp]
		\centering
		\caption{Impact evaluation of combining heuristic and exact methods for solving the pricing subproblem (best results in boldface) -- instances with $p_{\max}=1000$}
		\scriptsize
		\setlength{\tabcolsep}{0.8mm} 
		\begin{tabular}{lrHrrrrrrrrrrrrrrrr}
			\toprule
			\multirow{3}{*}{$n$} & \multirow{3}{*}{$m$} & \multirow{3}{*}{$\alpha$} &	& \multicolumn{7}{c}{LC$_A$} &       & \multicolumn{7}{c}{LC$_{2l}$}\\
			\cmidrule{5-11}\cmidrule{13-19}
			&	&	&	& \multicolumn{3}{c}{$H_\text{1-cycle}$} &       & \multicolumn{3}{c}{$H_{\text{elem}}$} &       & \multicolumn{3}{c}{$H_\text{1-cycle}$} &       & \multicolumn{3}{c}{$H_{\text{elem}}$} \\
			\cmidrule{5-7}\cmidrule{9-11}\cmidrule{13-15}\cmidrule{17-19}
			&  &  & & cols  & $\%$cols$_e$ & $t(s)$  &       & cols  & $\%$cols$_e$ & $t(s)$  &       & cols  & $\%$cols$_e$ & $t(s)$  &       & cols  & $\%$cols$_e$ & $t(s)$ \\
			\cmidrule{1-3}\cmidrule{5-7}\cmidrule{9-11}\cmidrule{13-15}\cmidrule{17-19}    {20} & 2     &       &       & 265.2 & 14.7  & 0.2   &       & 218.9 & 21.2  & 0.3   &       & 229.5 & \textbf{2.0} & \textbf{0.1} &       & \textbf{188.1} & 8.3   & 0.3 \\
			& 3     &       &       & 180.2 & 8.0   & 0.1   &       & 145.2 & 11.7  & 0.2   &       & 165.7 & \textbf{0.7} & \textbf{0.1} &       & \textbf{132.1} & 3.4   & 0.2 \\
			& 5     &       &       & 122.8 & 3.7   & 0.0   &       & 97.6  & 6.2   & 0.1   &       & 118.9 & \textbf{0.9} & \textbf{0.0} &       & \textbf{93.1} & 1.9   & 0.1 \\
			& 10    &       &       & 71.0  & 2.8   & \textbf{0.0} &       & 58.9  & 2.9   & 0.0   &       & 70.1  & \textbf{1.5} & 0.0   &       & \textbf{58.2} & 1.7   & 0.0 \\
			\cmidrule{1-3}\cmidrule{5-7}\cmidrule{9-11}\cmidrule{13-15}\cmidrule{17-19}    {50} & 2     &       &       & 2093.8 & 25.6  & 7.1   &       & 1743.0 & 33.6  & 16.2  &       & 1691.9 & \textbf{8.2} & \textbf{6.5} &       & \textbf{1451.8} & 20.0  & 17.1 \\
			& 3     &       &       & 1241.2 & 14.7  & 3.1   &       & 1011.4 & 18.6  & 7.6   &       & 1095.0 & \textbf{3.4} & \textbf{2.9} &       & \textbf{900.0} & 8.7   & 7.8 \\
			& 5     &       &       & 734.4 & 6.5   & 1.1   &       & 580.5 & 8.0   & 2.9   &       & 691.3 & \textbf{0.8} & \textbf{1.1} &       & \textbf{549.5} & 3.0   & 2.9 \\
			& 10    &       &       & 423.5 & 2.7   & 0.4   &       & 346.2 & 3.2   & 1.0   &       & 413.2 & \textbf{0.3} & \textbf{0.4} &       & \textbf{338.0} & 0.9   & 1.0 \\
			\cmidrule{1-3}\cmidrule{5-7}\cmidrule{9-11}\cmidrule{13-15}\cmidrule{17-19}    {100} & 2     &       &       & 10732.0 & 40.8  & 243.2 &       & 10461.6 & 49.4  & 488.8 &       & \textbf{8684.4} & \textbf{26.8} & \textbf{242.8} &       & 9624.9 & 44.2  & 563.0 \\
			& 3     &       &       & 5523.0 & 24.0  & 66.0  &       & 4462.1 & 30.9  & 135.6 &       & 4522.2 & \textbf{7.3} & \textbf{57.4} &       & \textbf{3698.6} & 17.0  & 134.7 \\
			& 5     &       &       & 2832.7 & 11.3  & 20.3  &       & 2226.4 & 14.7  & 45.5  &       & 2583.0 & \textbf{2.8} & \textbf{19.0} &       & \textbf{2017.1} & 5.8   & 45.4 \\
			& 10    &       &       & 1470.0 & 3.6   & 5.3   &       & 1244.8 & 4.5   & 14.3  &       & 1422.0 & \textbf{0.4} & \textbf{5.1} &       & \textbf{1204.3} & 1.3   & 14.3 \\
			\cmidrule{1-3}\cmidrule{5-7}\cmidrule{9-11}\cmidrule{13-15}\cmidrule{17-19}    \multicolumn{3}{l}{Avg.} &       & 2140.8 & 13.2  & 28.9  &       & 1883.1 & 17.1  & 59.4  &       & 1807.3 & \textbf{4.6} & \textbf{27.9} &       & \textbf{1688.0} & 9.7   & 65.6 \\
			\bottomrule
		\end{tabular}%
		\label{tab:CG_heur+exact_pmax1000}%
	\end{table}%
	It can be seen that, differently from the results for the instances with $p_{\max} = 100$, the use of LC$_{2l}$ combined with $H_\text{1-cycle}$ as heuristic method turns out to be the best alternative with regard to execution time. This result can be explained by the fact that $H_\text{1-cycle}$ runs in $\mathcal{O}{(nT)}$ time while $H_{\text{elem}}$ needs $\mathcal{O}{(n^2T)}$ time. However, it can also be observed that by using $H_{\text{elem}}$ as heuristic approach the number of required columns to prove the optimality is smaller than when $H_\text{1-cycle}$ is used. Hence, there exists a trade-off between execution times and number of generated columns. 
	
	{In the next experiments, reported in Tables \ref{tab:BP_MILP_20-100_pmax1000} and \ref{tab:BP_MILP_125-200_pmax1000}, we thus decided to test both B\&P methods that use LC$_{2l}$ with $H_\text{elem}$ and LC$_{2l}$ with $H_\text{1-cycle}$. The two B\&P methods are labeled as B\&P$_h$ for $H_{\text{elem}}$ and B\&P$_1$ for $H_\text{1-cycle}$. For the instances with up to $100$ jobs, it can be observed that  TI  is not able to manage instances with more than $50$ jobs due to memory limitations, while AF performs better, being able to tackle instances with up to $100$ jobs, especially when $m$ is large. The B\&P methods in turn, solve instances with up to $100$ jobs, but with average execution times slightly better than those of AF.}
	
	\begin{table}[htb]
		\centering
		\caption{Comparison of B\&P$_h$ and B\&P$_1$ performance with MILP models TI and AF - instances with $p_{\max}=1000$ and $n \in \{20,50,100\}$ (best results in boldface)}
		\scriptsize
		\setlength{\tabcolsep}{0.45mm} 
		\scalebox{0.9}{
			\begin{tabular}{lrrrrrrrrrrrrrrrrrrrrrrrrr}
				\toprule
				\multirow{3}{*}{$n$} & \multirow{3}{*}{$m$} & &  \multicolumn{5}{c}{TI} & & \multicolumn{5}{c}{AF} & & \multicolumn{5}{c}{B\&P$_h$} & & \multicolumn{5}{c}{B\&P$_1$}\\
				\cmidrule{4-8}\cmidrule{10-14}\cmidrule{16-20}\cmidrule{22-26}
				&  &  & \multirow{2}{*}{$opt$} & {$gap_{lp}$} & {$gap$} & \multirow{2}{*}{$nd$}  & \multirow{2}{*}{$t(s)$} &       & \multirow{2}{*}{$opt$} & {$gap_{lp}$} & {$gap$} & \multirow{2}{*}{$nd$}  & \multirow{2}{*}{$t(s)$} & & \multirow{2}{*}{$opt$} & {$gap_{lp}$} & {$gap$} & \multirow{2}{*}{$nd$}  & \multirow{2}{*}{$t(s)$} & & \multirow{2}{*}{$opt$} & {$gap_{lp}$} & {$gap$} & \multirow{2}{*}{$nd$}  & \multirow{2}{*}{$t(s)$} \\
				&  &  &  & {$(\%)$} & {$(\%)$} &  &  &       &   & {$(\%)$} & {$(\%)$} &  &  &       &   & {$(\%)$} & {$(\%)$} &  & & &   & {$(\%)$} & {$(\%)$} &  &   \\		
				\cmidrule{1-2}\cmidrule{4-8}\cmidrule{10-14}\cmidrule{16-20}\cmidrule{22-26}    {20} & 2     &       & \textbf{6} & 0.328 & \textbf{0.000} & 0.7   & 617.6 &       & \textbf{6} & 0.328 & \textbf{0.000} & \textbf{0.5} & 13.9  &       & \textbf{6} & \textbf{0.091} & \textbf{0.000} & 3.7   & 0.7   &       & \textbf{6} & \textbf{0.091} & \textbf{0.000} & 3.3   & \textbf{0.3} \\
				& 3     &       & \textbf{6} & \textbf{0.000} & \textbf{0.000} & \textbf{0.0} & 280.5 &       & \textbf{6} & \textbf{0.000} & \textbf{0.000} & \textbf{0.0} & 5.4   &       & \textbf{6} & \textbf{0.000} & \textbf{0.000} & \textbf{0.0} & 0.2   &       & \textbf{6} & \textbf{0.000} & \textbf{0.000} & \textbf{0.0} & \textbf{$<$0.1} \\
				& 5     &       & \textbf{6} & \textbf{0.000} & \textbf{0.000} & \textbf{0.0} & 148.6 &       & \textbf{6} & \textbf{0.000} & \textbf{0.000} & \textbf{0.0} & 2.7   &       & \textbf{6} & \textbf{0.000} & \textbf{0.000} & \textbf{0.0} & \textbf{$<$0.1} &       & \textbf{6} & \textbf{0.000} & \textbf{0.000} & \textbf{0.0} & \textbf{$<$0.1} \\
				& 10    &       & \textbf{6} & \textbf{0.000} & \textbf{0.000} & \textbf{0.0} & 67.2  &       & \textbf{6} & \textbf{0.000} & \textbf{0.000} & \textbf{0.0} & 1.1   &       & \textbf{6} & \textbf{0.000} & \textbf{0.000} & \textbf{0.0} & \textbf{$<$0.1} &       & \textbf{6} & \textbf{0.000} & \textbf{0.000} & \textbf{0.0} & \textbf{$<$0.1} \\
				\cmidrule{1-2}\cmidrule{4-8}\cmidrule{10-14}\cmidrule{16-20}\cmidrule{22-26}    {50} & 2     &       & 0     &       &       &       & m.lim &       & 3     & 0.122 & 0.057 & 56.5  & 946.9 &       & \textbf{6} & \textbf{0.079} & \textbf{0.000} & 8.5   & 81.5  &       & \textbf{6} & \textbf{0.079} & \textbf{0.000} & \textbf{8.3} & \textbf{37.5} \\
				& 3     &       & 0     &       &       &       & m.lim &       & \textbf{6} & 0.074 & \textbf{0.000} & 40.3  & 179.2 &       & 5     & \textbf{0.021} & 0.006 & 50.2  & 317.8 &       & \textbf{6} & \textbf{0.021} & \textbf{0.000} & \textbf{24.3} & \textbf{60.2} \\
				& 5     &       & 4     & 0.044 & 0.030 & \textbf{0.8} & 1515.1 &       & \textbf{6} & 0.044 & \textbf{0.000} & 441.2 & 238.6 &       & \textbf{6} & \textbf{0.026} & \textbf{0.000} & 33.5  & 18.3  &       & \textbf{6} & \textbf{0.026} & \textbf{0.000} & 29.0  & \textbf{8.3} \\
				& 10    &       & \textbf{6} & 0.002 & \textbf{0.000} & \textbf{0.3} & 531.5 &       & \textbf{6} & 0.002 & \textbf{0.000} & \textbf{0.3} & 7.7   &       & \textbf{6} & \textbf{0.001} & \textbf{0.000} & 2.7   & 1.3   &       & \textbf{6} & \textbf{0.001} & \textbf{0.000} & 4.0   & \textbf{0.5} \\
				\cmidrule{1-2}\cmidrule{4-8}\cmidrule{10-14}\cmidrule{16-20}\cmidrule{22-26}    {100} & 2     &       &       &       &       &       & m.lim       &       & \textbf{1} & 0.167 & 0.161 & 9.2   & \textbf{1635.2} &       & 0     & \textbf{0.142} & 0.132 & \textbf{9.0} & t.lim  &       & 0     & \textbf{0.142} & \textbf{0.131} & 12.2  & t.lim \\
				& 3     &       &       &       &       &       &   m.lim    &       & \textbf{1} & 0.102 & 0.100 & \textbf{13.7} & \textbf{1556.8} &       & 0     & \textbf{0.093} & 0.064 & 33.3  & t.lim  &       & 0     & \textbf{0.093} & \textbf{0.062} & 45.0  & t.lim \\
				& 5     &       &       &       &       &       &   m.lim    &       & 1     & 0.095 & 0.073 & \textbf{72.0} & 1523.5 &       & 1     & \textbf{0.079} & \textbf{0.056} & 194.2 & 1558.8 &       & \textbf{2} & \textbf{0.079} & \textbf{0.056} & 227.2 & \textbf{1314.5} \\
				& 10    &       &       &       &       &       &    m.lim   &       & \textbf{5} & 0.012 & \textbf{$<$0.001} & \textbf{163.2} & \textbf{637.0} &       & 4     & \textbf{0.007} & 0.001 & 748.2 & 719.0 &       & 4     & \textbf{0.007} & 0.001 & 820.7 & 643.5 \\
				\cmidrule{1-2}\cmidrule{4-8}\cmidrule{10-14}\cmidrule{16-20}\cmidrule{22-26}
				\multicolumn{2}{l}{Sum/Avg.$^*$} &       & 34    & 0.062 & 0.005 & 0.3   & 526.8 &       & 36    & 0.062 & \textbf{0.000} & 73.7 & 44.9 &       & \textbf{36}    & \textbf{0.020} & \textbf{0.000} & 6.6  & 3.4 &       & \textbf{36} & \textbf{0.020} & \textbf{0.000} & \textbf{6.1}  & \textbf{1.5} \\
				\multicolumn{2}{l}{Sum/Avg.} &       &     &  &  &    &  &       & 53    & 0.079 & 0.033 & \textbf{66.4} & 562.3 &       & 52    & \textbf{0.045} & 0.022 & 78.8  & 524.8 &       & \textbf{54} & \textbf{0.045} & \textbf{0.021} & 97.8  & \textbf{472.1} \\
				\bottomrule
				\multicolumn{26}{l}{$^*$Sum and average results for all instances, but the ones with $(n,m) \in \{(50,2),(50,3), (100,2), (100,3), (100,5), (100,10)\}$}
			\end{tabular}%
		}
		\label{tab:BP_MILP_20-100_pmax1000}%
	\end{table}%
	
	\begin{table}[htb]
		\centering
		\caption{Comparison of B\&P$_h$ and B\&P$_1$ performance with AF model - instances with $p_{\max}=1000$ and $n \in \{125,150\}$ (best results in boldface)}
		\scriptsize
		\setlength{\tabcolsep}{0.6mm} 
		\scalebox{0.9}{
			\begin{tabular}{lrrHHHHHHrrrrrrrrrrrrrrrrr}
				\toprule
				\multirow{3}{*}{$n$} & \multirow{3}{*}{$m$} & &  \multicolumn{5}{c}{} & & \multicolumn{5}{c}{AF} & & \multicolumn{5}{c}{B\&P$_h$} & & \multicolumn{5}{c}{B\&P$_1$}\\
				\cmidrule{4-8}\cmidrule{10-14}\cmidrule{16-20}\cmidrule{22-26}
				&  &  & \multirow{2}{*}{$opt$} & {$gap_{lp}$} & {$gap$} & \multirow{2}{*}{$nd$}  & \multirow{2}{*}{$t(s)$} &       & \multirow{2}{*}{$opt$} & {$gap_{lp}$} & {$gap$} & \multirow{2}{*}{$nd$}  & \multirow{2}{*}{$t(s)$} & & \multirow{2}{*}{$opt$} & {$gap_{lp}$} & {$gap$} & \multirow{2}{*}{$nd$}  & \multirow{2}{*}{$t(s)$} & & \multirow{2}{*}{$opt$} & {$gap_{lp}$} & {$gap$} & \multirow{2}{*}{$nd$}  & \multirow{2}{*}{$t(s)$} \\
				&  &  &  & {$(\%)$} & {$(\%)$} &  &  &       &   & {$(\%)$} & {$(\%)$} &  &  &       &   & {$(\%)$} & {$(\%)$} &  & & &   & {$(\%)$} & {$(\%)$} &  &   \\					
				\cmidrule{1-2}\cmidrule{4-8}\cmidrule{10-14}\cmidrule{16-20}\cmidrule{22-26}    {125} & 2     &       & 0     &       &       &       & m.lim &       & 0     & 16.885 & 16.884 & 4.0   & t.lim  &       & 0     & \textbf{0.209} & 0.209 & \textbf{1.2} & t.lim  &       & 0     & \textbf{0.209} & \textbf{0.207} & 2.2   & t.lim \\
				& 3     &       & 0     &       &       &       & m.lim &       & 0     & 0.176 & 0.173 & \textbf{7.3} & t.lim  &       & 0     & \textbf{0.157} & 0.141 & 14.3  & t.lim  &       & 0     & \textbf{0.157} & \textbf{0.139} & 21.8  & t.lim \\
				& 5     &       & 0     &       &       &       & m.lim &       & \textbf{1} & 0.141 & 0.140 & \textbf{11.3} & \textbf{1594.1} &       & \textbf{1} & \textbf{0.138} & \textbf{0.128} & 78.3  & 1634.6 &       & 0     & \textbf{0.138} & \textbf{0.128} & 114.3 & t.lim \\
				& 10    &       & 0     &       &       &       & m.lim &       & \textbf{2} & 0.080 & 0.078 & \textbf{68.2} & \textbf{1304.6} &       & \textbf{2} & \textbf{0.072} & \textbf{0.045} & 583.5 & 1388.8 &       & \textbf{2} & \textbf{0.072} & \textbf{0.045} & 793.5 & 1365.9 \\
				\cmidrule{1-2}\cmidrule{4-8}\cmidrule{10-14}\cmidrule{16-20}\cmidrule{22-26}    {150} & 2     &       & 0     &       &       &       & m.lim &       & 0     &       &       &       & m.lim &       & 0     & 0.282 & 0.282 & \textbf{0.2} & t.lim  &       & 0     & \textbf{0.281} & \textbf{0.253} & 0.7   & t.lim \\
				& 3     &       & 0     &       &       &       & m.lim &       & 0     &       &       &       & m.lim &       & 0     & 0.185 & 0.184 & \textbf{4.7} & t.lim  &       & 0     & \textbf{0.185} & \textbf{0.180} & 8.3   & t.lim \\
				& 5     &       & 0     &       &       &       & m.lim &       & 0     &       &       &       & m.lim &       & \textbf{1} & \textbf{0.224} & 0.212 & \textbf{43.3} & \textbf{1777.5} &       & 0     & \textbf{0.224} & \textbf{0.211} & 74.5  & t.lim \\
				& 10    &       & 0     &       &       &       & m.lim &       & 0     &       &       &       & m.lim &       & 0     & \textbf{0.097} & \textbf{0.094} & \textbf{547.0} & t.lim  &       & 0     & \textbf{0.097} & \textbf{0.094} & 681.7 & t.lim \\
				\cmidrule{1-2}\cmidrule{4-8}\cmidrule{10-14}\cmidrule{16-20}\cmidrule{22-26}    {200} & 2     &       & 0     &       &       &       & m.lim &       & 0     &       &       &       & m.lim &       & 0     & 0.405 & 0.405 & \textbf{0.0} & t.lim  &       & 0     & \textbf{0.341} & \textbf{0.341} & \textbf{0.0} & t.lim \\
				& 3     &       & 0     &       &       &       & m.lim &       & 0     &       &       &       & m.lim &       & 0     & 0.205 & 0.205 & \textbf{0.0} & t.lim  &       & 0     & \textbf{0.199} & \textbf{0.199} & 0.2   & t.lim \\
				& 5     &       & 0     &       &       &       & m.lim &       & 0     &       &       &       & m.lim &       & 0     & \textbf{0.272} & 0.267 & \textbf{9.5} & t.lim  &       & 0     & \textbf{0.272} & \textbf{0.266} & 19.2  & t.lim \\
				& 10    &       & 0     &       &       &       & m.lim &       & 0     &       &       &       & m.lim &       & 0     & \textbf{0.211} & \textbf{0.207} & \textbf{152.2} & t.lim  &       & 0     & \textbf{0.211} & \textbf{0.207} & 215.2 & t.lim \\
				\cmidrule{1-2}\cmidrule{4-8}\cmidrule{10-14}\cmidrule{16-20}\cmidrule{22-26}
				\multicolumn{2}{l}{Sum/Avg.$^*$} &       & 0     &       &       &       & m.lim &       & 3     & 4.321 & 4.318 & \textbf{23.5}  & \textbf{1624.7} &       & \textbf{3} & \textbf{0.144} & 0.131 & 169.3 & 1655.8 &       & 2     & \textbf{0.144} & \textbf{0.130} & 233.0 & 1691.5 \\
				\multicolumn{2}{l}{Sum/Avg.} &       &      &       &       &       &  &       &      &  &  &   &  &       & \textbf{4} & 0.205 & 0.198 & \textbf{119.5} & \textbf{1750.1} &       & 2     & \textbf{0.199} & \textbf{0.189} & 161.0 & 1763.8 \\
				\bottomrule
				\multicolumn{23}{l}{$^*$Sum and average results for the instances with $n = 125$}
			\end{tabular}%
		}
		\label{tab:BP_MILP_125-200_pmax1000}%
	\end{table}%
	
	Concerning the experiments on the larger instances, {TI reached memory limit on all of them and was consequently disregarded from Table \ref{tab:BP_MILP_125-200_pmax1000}. From this table,} it can be noticed that for instances with $125$ jobs, AF, B\&P$_h$ and B\&P$_1$ have a similar performance. However, AF was not able to tackle the instances with $150$ jobs due to memory limit. In these cases, the B\&P methods are a more robust alternative, being able to prove the optimality of one instance and to provide very low gaps for the other ones, {always below 0.4\% on each group of instances}. In particular, B\&P$_h$ performs better than  B\&P$_1$ on the larger instances, solving $4$ instances to the proven optimality instead of 2. 
	
	\section{Concluding remarks} \label{sec:coclusion_PrwC}
	
	{In this work, we investigated the problem of scheduling a set of jobs that are released over the time on a number of identical parallel machines, with the aim of minimizing the total weighted completion time. This problem, denoted as {\PrwC}, is of high interest because can model several applications in practice. We proposed two MILP formulations, a time-indexed and an arc-flow, and a tailored B\&P algorithm. Our experiments showed that for small and medium-sized instances all attempted methods performed quite well, with the arc-flow model achieving the best performance in terms of number of proven optimal solutions. {Regarding the large-sized instances involving more than $200$ jobs, or instances with large processing times, we notice that the MILP formulations suffers from memory limit due to their pseudo-polynomial sizes, whereas the B\&P remains a robust and interesting alternative being able to provide very low gaps in a reasonable amount of time. Using the proposed methods we were able to solve for the first time instances with up to $200$ jobs to the proven optimality.
			
			Concerning  future research direction we observe that  a large number of instances, especially involving large processing times, remain unsolved, thus  leaving space for further investigation. Moreover, we envisage the need of a deeper study of the proposed models aiming to reduce  their sizes by considering combinatorial properties of the problem, as well as the use of advanced concepts for the design of the B\&P algorithm, such as, the introduction of cutting planes, the use of stabilization techniques, the implementation of an early termination approach for the CG algorithm and the use of a strong branching scheme. The developed techniques could be also adapted to tackle other problem variants, especially those involving due dates and the minimization of (weighted) tardiness, lateness or just-in-time objective functions.
		}

		
		
		\section*{Acknowledgments}
		This research was partially funded by the CNPq - Conselho Nacional de Desenvolvimento Cient{\'i}fico e Tecnol{\'o}gico, Brazil, grant No. 234814/2014-4 and by University of Modena and Reggio Emilia, under grant FAR 2018 Analysis and optimization of health care and pharmaceutical logistic processes.
		
\bibliographystyle{mmsbib}
\bibliography{ref}

\end{document}